\documentclass[11pt,reqno]{amsart}
\usepackage{bbm}
\usepackage{amsfonts}
\usepackage{cases}
\usepackage{mathrsfs}
\usepackage{amsmath}
\usepackage{amssymb, amsmath}
\usepackage{graphicx}

\newcommand{\newcom}{\newcommand}

\newcom{\al}{\alpha}
\newcom{\be}{\beta}
\newcom{\eps}{\epsilon}
\newcom{\ga}{\gamma}
\newcom{\Ga}{\Gamma}
\newcom{\ka}{\kappa}
\newcom{\Lam}{\Lambda}
\newcom{\lam}{\lambda}
\newcom{\la}{\lambda}
\newcom{\Om}{\Omega}
\newcom{\om}{\omega}
\newcom{\Si}{\Sigma}
\newcom{\si}{\sigma}
\newcom{\tht}{\theta}
\newcom{\dtri}{\nabla}
\newcom{\tri}{\triangle}
\newcom{\oo}{\infty}
\newcom{\vphi}{\varphi}
\newcom{\cA}{{\mathcal A}}
\newcom{\cB}{{\mathcal B}}
\newcom{\cC}{{\mathcal C}}
\newcom{\cD}{{\mathcal D}}
\newcom{\cE}{{\mathcal E}}
\newcom{\cF}{{\mathcal F}}
\newcom{\cG}{{\mathcal G}}
\newcom{\cL}{{\mathcal L}}
\newcom{\cM}{{\mathcal M}}
\newcom{\cP}{{\mathcal P}}
\newcom{\cS}{{\mathcal S}}
\newcom{\cQ}{{\mathcal Q}}
\newcom{\caly}{{\mathcal Y}}
\newcom{\calZ}{{\mathcal Z}}
\newcom{\bfz}{{\bf Z}}
\newcom{\R}{\Bbb R}
\newcom{\N}{\Bbb N}
\newcom{\Z}{\Bbb Z}
\newcom{\C}{\Bbb C}
\newcom{\E}{\Bbb E}

\def\dv{\mbox{div}}

\newcom{\hn}{{\bf H}^n}
\newcom{\hnn}{{\mathbf H}^{n'}}
\newcom{\ulzs}{u^\lam_{z,s}}
\newcom{\Hl}{{{\cal  H}_\lam}}
\newcom{\fal}{F_{\al, \lam}}
\newcom{\Dh}{\Delta_{{\mathbf H}^n}}
\newcom{\fgl}{F_{\g, \lam}}

\newcom{\f}{\frac}
\newcom{\di}{\displaystyle\int}
\newcom{\ds}{\displaystyle\sum}
\newcom{\dl}{\displaystyle\lim}
\newcom{\ov}{\overline}
\newcom{\sset}{\subset}
\newcom{\wt}{\widetilde}
\newcom{\pa}{\partial}
\newcom{\p}{\partial}
\newcom\na{\nabla}
\newcom{\co}{\cdot}
\newcom{\suml}{\sum\limits}
\newcom{\supl}{\sup\limits}
\newcom{\intl}{\int\limits}
\newcom{\infl}{\inf\limits}
\newcom{\disp}{\displaystyle}
\newcom{\non}{\nonumber}
\newcom{\no}{\noindent}
\newcom{\QED}{$\square$}
\newcom{\tphi}{\tilde\phi}
\newcom{\curl}{\rm{curl}}
\newcom{\qqquad}{\qquad\qquad\qquad}

\def\ef{\hphantom{MM}\hfill\llap{$\square$}\goodbreak}

\def\eqdefa{\buildrel\hbox{\footnotesize def}\over {=\!\!=}}

\newtheorem{athm}{\bf \t}[section]
\newenvironment{thm} [1] {\def\t{#1}\begin{athm} \bf \rm} {\end{athm}}
\newcom{\bthm}{\begin{thm}}\newcom{\ethm}{\end{thm}}

\newcom{\beq}{\begin{equation}}
\newcom{\eeq}{\end{equation}}
\newcom{\ben}{\begin{eqnarray}}
\newcom{\een}{\end{eqnarray}}
\newcom{\beno}{\begin{eqnarray*}}
\newcom{\eeno}{\end{eqnarray*}}
\newcom{\bali}{\begin{aligned}}
\newcom{\eali}{\end{aligned}}

\newtheorem{Theorem}{Theorem}[section]
\newtheorem{Definition}[Theorem]{Definition}
\newtheorem{Proposition}[Theorem]{Proposition}
\newtheorem{Lemma}[Theorem]{Lemma}

\newtheorem{Remark}[Theorem]{Remark}

\numberwithin{equation}{section}

\begin{document}

\title[Ill-posedness of the compressible Navier-Stokes equations]
{\small On the ill-posedness of the compressible Navier-Stokes equations in the critical Besov spaces}

\author{Qionglei Chen}
\address{Institute of Applied Physics and Computational Mathematics,P.O. Box 8009, Beijing 100088, P. R. China}
\email{chen\_qionglei@iapcm.ac.cn}

\author{Changxing Miao}
\address{Institute of Applied Physics and Computational Mathematics,P.O. Box 8009, Beijing 100088, P. R. China}
\email{miao\_changxing@iapcm.ac.cn}

\author{Zhifei Zhang}
\address{School of Mathematical Sciences, Peking University, 100871, P. R. China}
\email{zfzhang@math.pku.edu.cn}


\date{\today}



\begin{abstract}
We prove the ill-posedness of the 3-D  baratropic  Navier-Stokes
equation for the initial density and velocity belonging to the
critical Besov space $(\dot{B}^{\f
3p}_{p,1}+\bar{\rho},\,\dot{B}^{\f 3p-1}_{p,1})$ for $p>6$ in the
sense that a ``norm inflation" happens in finite time, here
$\bar{\rho}$ is a positive constant. Our argument also shows that
the compressible viscous heat-conductive flows is ill-posed for the
initial density, velocity and temperature belonging to the critical
Besov space $(\dot{B}^{\f 3p}_{p,1}+\bar{\rho},\,\dot{B}^{\f
3p-1}_{p,1},\,\dot{B}^{\f 3p-2}_{p,1})$ for $p>3$. These results
shows that the compressible Navier-Stokes equations are ill-posed in
the smaller critical spaces compared with the incompressible
Navier-Stokes equations.
\end{abstract}

\maketitle

\section{Introduction}
The full compressible Navier-Stokes equations read as follows
\begin{align}\label{equ:CCNS}
\left\{
\begin{aligned}
&\p_t\rho+\textrm{div}(\rho u)=0,\\
&\p_t(\rho u)+\textrm{div}(\rho u\otimes u)=\textrm{div}\tau, \\
&\p_t\big(\rho u(e+\frac{|u|^2}{2})\big)+\textrm{div}\big(\rho
u(e+\frac{|u|^2}{2})\big)=\textrm{div}(\tau\cdot u+\kappa\nabla\theta),
\end{aligned}
\right.
\end{align}
where $\rho(t,x)$, $u(t,x)$, $e(t,x)$ denote the density,
velocity  of the fluid and the internal energy per unit mass respectively,  $\kappa>0$ is the thermal conduction parameter, and $\theta$ is the temperature.
The internal stress tensor $\tau$ is given by
$$\tau=2\nu {\rm D}(u)+(\lambda\textrm{div}u-p)I,$$
where ${\rm D}(u)=\frac12(\nabla u+\nabla u^{\!\top})$, the constants $\mu,\lambda$ are the viscosity coefficients satisfying
\begin{align}
\mu>0\quad \textrm{and} \quad \lambda+2\mu>0.\nonumber
\end{align}

For the ideal gas, $e=c_V\theta,\, P=\rho R\theta$ for some constants $c_V>0, R>0$. In such case, the system \eqref{equ:CCNS}
can be rewritten as
\begin{align}\label{equ:cNS}
\left\{
\begin{aligned}
&\p_t\rho+\textrm{div}(\rho u)=0,\\
&\p_t(\rho u)+\textrm{div}(\rho u\otimes u)-\mu\Delta u-(\lambda+\mu)\na\textrm{div}u+\na P=0, \\
&c_V(\p_t(\rho \theta)+\textrm{div}(\rho
u\theta))-\kappa\Delta\theta+P\textrm{div}u=\frac{\mu}{2}|\nabla
u+(\nabla u)^{\!\top}|^2+
\lambda|\textrm{div}u|^2.
\end{aligned}
\right.
\end{align}
Here we denote by $|A|^2$ the trace of the matrix $AA^{\!\top}$.

When the pressure depends  only on the density, we get the
baratropic Navier-Stokes equations
\begin{align}\label{equ:bNS}
\left\{
\begin{aligned}
&\p_t\rho+\textrm{div}(\rho u)=0,\\
&\p_t(\rho u)+\textrm{div}(\rho u\otimes u)-\mu\Delta u-(\lambda+\mu)\na\textrm{div}u+\na P=0.
\end{aligned}
\right.
\end{align}

In this paper, we are concerned with the Cauchy problem of the
system (\ref{equ:cNS}) and (\ref{equ:bNS}) in $\R^+\times \R^3$
together with the initial data \beno (\rho,u,\theta)|_{t=0}=(\rho_0,
u_0, \theta_0),\quad\mbox{and}\quad (\rho,u)|_{t=0}=(\rho_0, u_0),
\eeno respectively.

The local existence and uniqueness of smooth solutions for the system
(\ref{equ:cNS}) were proved by Nash \cite{Nash} for smooth initial
data without vacuum. Matsumura-Nishida\cite{Mat-Nish} proved that the solution is global in time for the data close to equilibrium.
For the general case, the question of whether smooth solutions blow up in finite time is widely open, even in two dimensional case.
For the initial density with compact support, Xin \cite{Xin} proved that any non-zero smooth solutions of (\ref{equ:cNS}) will blow up in finite time.
Recently, Sun-Wang-Zhang \cite{SWZ-JMPA, SWZ} showed that smooth solution will not blow up as long as the upper bound of
the density( and temperature for (\ref{equ:cNS})) is bounded. We refer to the seminal books
\cite{Lions, Feir1} and references therein for the global existence of weak solutions.

Motivated by Fujita-Kato's result on the incompressible
Navier-Stokes equations \cite{Fujita}, Danchin applied Fourier
analysis method to study the well-posedness for the compressible
Navier-Stokes equations in critical spaces. Let us make it precise.
It is easy to check that if $(\rho,u,\theta)$ is a solution of
(\ref{equ:cNS}), then
\begin{align*}
(\rho_\lambda(t,x), u_\lambda(t,x),\theta_\lambda(t,x))\eqdefa
\bigl(\rho(\lambda^2 t,\lambda x),\lambda u(\lambda^2 t,\lambda x),\lambda^2 \theta(\lambda^2 t,\lambda x)\bigr),\quad\lambda>0
\end{align*}
is also a solution of (\ref{equ:cNS})  provided the pressure law has been changed into $\lambda^2 P$.
A functional space is called critical if the associated norm is invariant under the transformation
$
(\rho,u,\theta)\longrightarrow (\rho_\lambda,u_\lambda,\theta_\lambda)
$(up to a constant independent of $\lambda$). Roughly speaking, the system (\ref{equ:cNS}) is locally well-posed for the initial data
\begin{align*}
(\rho_0-\bar{\rho},u_0,\theta_0)\in\dot B^{\f 3 p}_{p,1}\times \bigl(\dot
B^{\f 3 p-1}_{p,1}\bigr)^3\times\dot B^{\f 3 p-2}_{p,1}\quad \textrm{with  } p<3;
\end{align*}
and the system (\ref{equ:bNS}) is locally well-posed for the initial
data
\begin{align*}
(\rho_0-\bar{\rho},u_0)\in\dot B^{\f 3 p}_{p,1}\times \bigl(\dot
B^{\f 3 p-1}_{p,1}\bigr)^3\quad \textrm{with  } p<6.
\end{align*}
Here $\dot B^s_{p,q}$ is the homogeneous Besov space, see Definition
\ref{Def:Bessov}. Furthermore, the system is globally well-posed if
the initial data is small in the critical Besov space with $p=2$. We
refer to \cite{Danchin-inven, Danchin-arma, Danchin-cpde01,
Danchin-cpde07, CMZ-Revista} and references therein. Recently,
Chen-Miao-Zhang and Charve-Danchin \cite{CMZcpam,Charve-Danchin}
proved that the system (\ref{equ:bNS}) is  globally well-posed for
the small initial data in the hybrid critical Besov spaces, in which
the part of high frequency of the initial data lies in the critical
Besov spaces with $p>3$, and  the part of low frequency lies in the
critical Besov spaces with $p=2$. This allows us to choose the
highly oscillating initial velocity like
$\sin(x_1/\varepsilon)\varphi(x)$ since it is small in the hybrid
critical Besov space.

A natural question is  whether the system (\ref{equ:cNS}) and (\ref{equ:bNS}) is well-posed in the critical Besov spaces with $p\ge 3$ and $p\ge 6$
respectively. Recently, for the incompressible Navier-Stokes equations, Bourgain-Pavlovi\'{c}\cite{Bourgain} and Germain \cite{Germain} proved the ill-posedness in  the largest critical space $\dot B^{-1}_{\infty,\infty}$. Motivated by \cite{Bourgain},
we prove that the system \eqref{equ:bNS} is ill-posed in the critical Besov spaces with $p>6$.

\bthm{Theorem}\label{thm:illposed-bNS}
Let $\bar{\rho}$  be a positive constant and $p>6$. For any $\delta>0$,
there exists  initial data $(\rho_0, u_0)$ satisfying
$$\|\rho_0-\bar{\rho}\|_{\dot{B}^{\frac3p}_{p,1}}\le \delta, \quad \|u_0\|_{\dot{B}^{\frac3p-1}_{p,1}}\le \delta$$
such that a solution $(\rho, u)$ to the system \eqref{equ:bNS} satisfies
$$\|u(t)\|_{\dot{B}^{\frac3p-1}_{p,1}}>\frac1\delta$$
for some $0<t<\delta$.
\ethm

While, the system \eqref{equ:cNS} is ill-posed in the critical Besov spaces with $p>3$.

\bthm{Theorem}\label{thm:illposed-heat}
Let $\bar{\rho}$  be a positive constant and $p>3$. For any $\delta>0$,
there exists initial data $(\rho_0, u_0, \theta_0)$ satisfying
$$\|\rho_0-\bar{\rho}\|_{\dot{B}^{\frac3p}_{p,1}}\le \delta, \quad \|u_0\|_{\dot{B}^{\frac3p-1}_{p,1}}\le \delta,\quad \|\theta_0\|_{\dot{B}^{\frac3p-2}_{p,1}}\le \delta$$
such that a solution $(\rho, u, \theta)$  to the system \eqref{equ:cNS} satisfies
$$\|\theta(t)\|_{\dot{B}^{\frac3p-2}_{p,1}}>\frac1\delta$$
for some $0<t<\delta$.
\ethm

\begin{Remark}For the
incompressible Navier-Stokes equation, the ill-posedness  is proved
in the largest  critical  space $\dot B^{-1}_{\infty,\infty}$ by
Bourgain-Pavlovi\'{c}\cite{Bourgain}. Generally speaking, the
ill-posedness is  easier to get if we work in the comparatively
larger space. However, for the more complex compressible
Navier-Stokes equations, noting that
$$\dot{B}^{\frac3p-1}_{p,1}\hookrightarrow
\dot{B}^{\frac3p-1}_{p,\infty}\hookrightarrow \dot
B^{-1}_{\infty,\infty}$$ for $p<\infty$, our results show that the
ill-posedness is established  even in the much better critical
spaces for $p>6$ or $p>3$ which have the same scaling with the
largest critical space.


\end{Remark}

\begin{Remark} Theorem  \ref{thm:illposed-heat} implies that it seems impossible to generate a global solution to \eqref{equ:cNS}
for the highly oscillating initial velocity as in \cite{Charve-Danchin, CMZcpam}. The question whether the system (\ref{equ:bNS})(or (\ref{equ:cNS}))
is well-posed in the critical Besov Space with $p=6$(or p=3) is still open.
\end{Remark}

For the baratropic Navier-Stokes equations, the mechanism leading to
the ill-posedness comes from the high-high frequency interaction of
the nonlinear term $u\cdot\na u$. For the viscous heat-conductive
flows, the mechanism leading to the ill-posedness comes from the
high-high frequency interaction of the strong nonlinear terms
$|\nabla u+(\nabla u)^{\!\top}|^2$ and $|\textrm{div}u|^2$ in the
temperature equation. Roughly speaking,  for the system
\eqref{equ:bNS}, we first decompose the velocity $u$ into \beno
u(t)=U_0(t)+U_1(t)+U_2(t), \eeno where $U_0(t)=e^{t\Delta}u_0$ and
$U_1(t)=\int_0^te^{(t-\tau)\Delta}U_0\cdot\na U_0(\tau)d\tau$.
Secondly, we construct a suitable combination of plane waves for the
initial velocity such that $\|U_1(t)\|_{\dot B^{\f3p-1}_{p,1}}$ is
big for some time $t>0$, while $\|U_1\|_{L^1_t\dot
B^{\frac3q+1}_{q,1}\cap \widetilde{L}^2_t\dot B^{\frac3q}_{q,1}}$ is
small for some $q<6$. Lastly, we will prove the remainder term
$\|U_2\|_{\widetilde{L}^\infty_T\dot B^{\frac3q-1}_{q,1}\cap
{L}^1_T\dot B^{\frac3q+1}_{q,1}}$ is also small by subtle estimates.

\section{Some tools of Littlewood-Paley analysis}

\subsection{Littlewood-Paley decomposition and Besov sapces}
Choose a radial function  $\varphi \in {\mathcal S}(\R^3)$ supported in
${\mathcal C}=\{\xi\in\R^3,\, \frac{3}{4}\le|\xi|\le\frac{8}{3}\}$ such
that
\begin{align*}
\sum_{j\in\Z}\varphi(2^{-j}\xi)=1 \quad \textrm{for
all}\,\,\xi\neq 0.
\end{align*}
The frequency localization operator
$\Delta_j$ and $S_j$ are defined by
\begin{align}
\Delta_jf=\varphi(2^{-j}D)f,\quad S_jf=\sum_{k\le
j-1}\Delta_kf\quad\mbox{for}\quad j\in \Z, \nonumber
\end{align}
respectively. The Fourier transform of $f$ is denoted by
$\widehat{f}$ or $\mathcal{F}f$, the inverse by $\mathcal{F}^{-1}f$.
We denote by ${\mathcal Z'}(\R^3)$  the dual space of ${\mathcal
Z}(\R^3)=\{f\in {\mathcal S}(\R^3);\,D^\alpha \widehat{f}(0)=0;
\forall\alpha\in\N^3 \,\mbox{multi-index}\}$, which  can also be
identified by the quotient space of ${\mathcal S'}(\R^3)/{\mathcal
P}$ with the polynomials space ${\mathcal P}$. Let us introduce the
homogeneous Besov space.

\begin{Definition}\label{Def:Bessov} Let $\sigma\in\R$, $1\le
r,s\le+\infty$. The homogeneous Besov space $\dot{B}^{\sigma}_{r,s}$
is defined by
$$\dot{B}^{\sigma}_{r,s}\eqdefa\{f\in {\mathcal Z'}(\R^3):\,\|f\|_{\dot{B}^{\sigma}_{r,s}}<+\infty\},$$
where \begin{align*}
\|f\|_{\dot{B}^{\sigma}_{r,s}}=\Bigl\|\big\{2^{k\sigma}
\|\Delta_kf\|_{L^r}\big\}_j\Bigr\|_{\ell^s}.\end{align*}
\end{Definition}

We next introduce the Chemin-Lerner type space
$\widetilde{L}^\rho_T\dot{B}^{\sigma}_{r,s}$.

\begin{Definition}Let $\sigma\in\R$, $1\le
\rho,r,s\le+\infty$, $0<T\le+\infty$. The space
$\widetilde{L}^\rho_T\dot{B}^\sigma_{r,s}$ is defined as the set of
all the distributions $f$ satisfying
$$\|f\|_{\widetilde{L}^\rho_T\dot{B}^{\sigma}_{r,s}}\eqdefa \Bigl\|\big\{2^{k\sigma}
\|\Delta_kf(t)\|_{L^\rho(0,T;L^r)}\big\}_k\Bigr\|_{\ell^s}<\infty.
$$
Obviously, $
\widetilde{L}^1_T(\dot{B}^\sigma_{r,1})=L^1_T(\dot{B}^\sigma_{r,1}).
$
\end{Definition}

Next we recall the estimates of the linear transport equation and
heat equation in Besov spaces which will be used in the subsequence.

\begin{Proposition}\label{Prop:transport}
Let $T>0$, $\sigma\in (-3\min(\frac1r,\frac1{r'}), 1+\frac 3r]$, and
$1\le r\le+\infty$. Let $v$ be a vector field so that $\nabla v\in
L^1_T(\dot{B}^{\frac{3}{r}}_{r,1})$. Assume that $u_0\in
\dot{B}^{\sigma}_{r,1},$ $f\in L^1_T(\dot{B}^{\sigma}_{r,1})$ and
$u$ is the solution of
\begin{align}
\bigg\{\begin{aligned}
&\partial_t u+v\cdot \nabla u =f,\\
&u(0,x)=u_0.\non
\end{aligned}
\bigg.\end{align} Then there holds for $t\in[0,T]$, \begin{align*}
\|u\|_{\widetilde{L}^\infty_t\dot{B}^{s}_{r,1}}\le e^{CV(t)}\Big(
\|u_0\|_{\dot{B}^{\sigma}_{r,1}}+\int_0^t
e^{-CV(\tau)}\|f(\tau)\|_{\dot{B}^{\sigma}_{r,1}}d\tau\Big),
\end{align*}
where $V(t)=\int_0^t\|\nabla
v(\tau)\|_{\dot{B}^{\frac{3}{r}}_{r,1}}d\tau.$
\end{Proposition}

\begin{Proposition}\label{Prop:heatflow}
Let $T>0$, $\sigma\in\R$ and $1\le r\le\infty$. Assume that
$u_0\in\dot B^\sigma_{r,1}$ and $f\in \widetilde{L}^\rho_T\dot
B^{\sigma-2+\frac2\rho}_{r,1}$. If $u$ is the solution of the heat
equation
\begin{align*}
\bigg\{\begin{aligned}
&\partial_t u-\mu\Delta u=f,\\
&u(0,x)=u_0,
\end{aligned}
\bigg.
\end{align*}
with $\mu>0$, then for all $\rho_1\in[\rho,\infty]$, it holds that
\begin{align}
\mu^{\frac1{\rho_1}}\|u\|_{\widetilde{L}^{\rho_1}_T\dot
B^{\sigma+\frac2{\rho_1}}_{r,1}}\le C\big(\|u_0\|_{\dot
B^{\sigma}_{r,1}}+\|f\|_{\widetilde{L}^{\rho}_T\dot
B^{\sigma-2+\frac2{\rho}}_{r,1}}\big).\non
\end{align}
\end{Proposition}

We refer to \cite{Bah-Che-Dan} for more details.

\subsection{Nonlinear estimates in Besov space}

Let us first recall some classical product estimates in Besov spaces from \cite{Bah-Che-Dan}.

\begin{Lemma}\label{Lem:binestis>0} Let $T>0$, $\sigma>0$ and $1\le r, \rho\le \infty$. Then it holds that
\begin{align*}
\|fg\|_{\widetilde{L}^\rho_T\dot B^{\sigma}_{r,1}} \le
C\big(\|f\|_{L^\infty_T(L^\infty)}\|g\|_{\widetilde{L}^\rho_T\dot
B^{\sigma}_{r,1}}+
\|g\|_{L^\infty_T(L^\infty)}\|f\|_{\widetilde{L}^\rho_T\dot
B^{\sigma}_{r,1}}\big).
\end{align*}
\end{Lemma}
\begin{Lemma}\label{Lem:binesti} Let $T>0$, $\sigma_1, \sigma_2\le \frac{3}{r},\, \sigma_1+\sigma_2>3\max (0,\frac2r-1)$ and $1\le r, \rho, \rho_1,\rho_2 \le \infty$
with $\f 1 \rho=\f 1{\rho_1}+\f1 {\rho_2}$. Then it holds that
\begin{align*} \|fg\|_{\widetilde{L}^\rho_T\dot B^{\sigma_1+\sigma_2-
\frac{3}{r}}_{r,1}} \le C\|f\|_{\widetilde{L}^{\rho_1}_T\dot
B^{\sigma_1}_{r,1}}\|g\|_{\widetilde{L}^{\rho_2}_T\dot
B^{\sigma_2}_{r,1}}.
\end{align*}
\end{Lemma}
\begin{Lemma}\label{Lem:nonesti}
Let $T>0$, $\sigma>0$ and $1\le r,\rho\le \infty$. Assume that $F\in
W^{[\sigma]+3,\infty}_{loc}(\R)$ with  $F(0)=0$. Then for any $f\in
L^\infty\cap \dot B^\sigma_{r,1}$, we have
\begin{align*}
\|F(f)\|_{\widetilde{L}^\rho_T\dot B^\sigma_{r,1}}\le
C\bigl(1+\|f\|_{L^\infty_T(L^\infty)}\bigr)^{[\sigma]+2}\|f\|_{\widetilde{L}^\rho_T\dot
B^\sigma_{r,1}}.
\end{align*}
\end{Lemma}

We will use  Bony's decomposition
\begin{align}\label{Bonydecom}
fg=T_{f}g+T_gf+R(f,g),
\end{align}
where
$$T_fg=\sum_{j\in\Z}S_{j-1}f\Delta_jg, \quad R(f,g)=\sum_{|j'-j|\le1}\Delta_jf{\Delta}_{j'}g.$$

We need the following estimates for the paraproduct $T_fg$ and $R(f,g)$.

\begin{Lemma}\label{lem:biest-n}
Let $\sigma, \al\in\R$ and $1\le r, b, \rho, \rho_1,
\rho_2\le\infty$ with $\f 1\rho=\f1 {\rho_1}+\f 1 {\rho_2}$. Then we
have

\begin{itemize}

\item[1.] if $\alpha\ge0$, then
\begin{align}
\|T_fg\|_{\widetilde{L}^\rho_T\dot B^{\sigma}_{b,1}}\le
\|f\|_{\widetilde{L}^{\rho_1}_T\dot B^{\f3r-\al}_{r,1}}
\|g\|_{\widetilde{L}^{\rho_2}_T\dot B^{\sigma+\al}_{b,1}};\non
\end{align}

\item[2.] if $\f3r-\f3b+\al\ge 0$ and $r\ge b$, then

\begin{align}
\|T_fg\|_{\widetilde{L}^\rho_T\dot B^{\sigma}_{b,1}}\le
C\|f\|_{\widetilde{L}^{\rho_1}_T\dot B^{\f3q-\al}_{b,1}}
\|g\|_{\widetilde{L}^{\rho_2}_T\dot B^{\sigma+\f 3r-\f3
b+\al}_{r,1}};\non
\end{align}

\item[3.] if $\frac3r+\sigma>0$, then
\begin{align}\label{equ:productIII}
\|R(f,g)\|_{\widetilde{L}^\rho_T\dot B^{\sigma}_{b,1}}\le
\|f\|_{\widetilde{L}^{\rho_1}_T\dot B^{\f3r-\al}_{r,1}}
\|g\|_{\widetilde{L}^{\rho_2}_T\dot B^{\sigma+\al}_{b,1}}.\non
\end{align}

\end{itemize}
\end{Lemma}

\no{\bf Proof.}\, Due to $\f3r-\f3b+\al\ge 0$, we infer from
H\"{o}lder's inequality that
\begin{align}
\|T_fg\|_{\widetilde{L}^\rho_T\dot B^{\sigma}_{b,1}} &\le \sum_{k\in
\mathbb{Z}}\|S_{k-1}f\|_{L^{\rho_2}_T(L^\frac{rb}{r-b})}\|\Delta_kg\|_{L^{\rho_1}_T(L^r)}2^{k\sigma}\nonumber\\&\le
\sum_{k'\le k-1}2^{(\frac3b-\al)
k'}\|\Delta_{k'}f\|_{L^{\rho_1}_T(L^b)}2^{(\frac3r-\frac3b+\al)k'}
\|\Delta_kg\|_{L^{\rho_2}_T(L^r)}2^{k\sigma}\nonumber\\& \le
C\|f\|_{\widetilde{L}^{\rho_1}_T\dot B^{\f3b-\al}_{b,1}}
\|g\|_{\widetilde{L}^{\rho_2}_T\dot B^{\sigma+\f 3r-\f3
b+\al}_{r,1}}.\non
\end{align}
Here we used the fact that \beno
\|\Delta_{k'}f\|_{L^{\frac{rb}{r-b}}}\le
C2^{\f3rk'}\|\Delta_{k'}f\|_{L^b}, \eeno which can be deduced from
Young's inequality. This proves (ii). The proof of (i) and (iii) is
similar.\ef

\section{Ill-posedness of the baratropic Navier-Stokes equations}

\subsection{Reformualtion of the equation}
We introduce the new unknowns
$$
a=\frac{\rho}{\bar\rho}-1,\quad h=\Lambda^{-1}\textrm{div}u,\quad
\Omega=\Lambda^{-1}{\curl}u,
$$
where $\Lambda^s f\eqdefa \mathcal{F}^{-1}(|\xi|^s\hat f(\xi))$. We
rewrite the velocity $u$ as follows
\begin{align*}
u=-\Lambda^{-1}\na h+\Lambda^{-1}{\curl}\Omega.
\end{align*}
The system \eqref{equ:bNS} can be rewritten as
\begin{equation}\label{equ:cNStrans}
\left\{
\begin{aligned}{}
&\p_ta+u\cdot\na a+\dv u(1+a)=0,\\
&\p_th-\bar{\nu}\Delta h=-\Lambda^{-1}\dv\big(u\cdot\na u+L(a)\mathcal{A}u+K(a)\nabla a\big),\\
&\p_t\Omega-\bar{\mu}\Delta \Omega=-\Lambda^{-1}{\curl}\big(u\cdot\na u+L(a)\mathcal{A}u\big),\\
&(a,h,\Omega)|_{t=0}=(a_0,h_0,\Omega_0),
\end{aligned}
\right.
\end{equation}where $$\mathcal
{A}=\bar{\mu}\Delta+(\bar{\lambda}+\bar{\mu})\na \dv, \quad
K(a)=\f{P'(\bar{\rho}(1+a))}{1+a}\quad\hbox{and}\quad
L(a)=\frac{a}{1+a}$$ with $\bar{\mu}=\f{\mu}{\bar{\rho}}$,
$\bar{\lambda}=\f{\lambda}{\bar{\rho}}$, and $\bar{\nu}=\bar{\lambda}+2\bar{\mu}$.
Hence, we obtain
\begin{align*}
&h(t,x)=e^{\bar\nu\Delta t}h_0-\int_0^te^{{\bar\nu(t-\tau)\Delta }}
\Lambda^{-1}\dv\,\big(u\cdot\na u+L(a)\mathcal{A}u+K(a)\nabla a\big)d\tau,\\
&\Omega(t,x)=e^{\bar\mu\Delta
t}\Omega_0-\int_0^te^{{\bar\mu(t-\tau)\Delta }}
\Lambda^{-1}{\curl}\,\big(u\cdot\na u+L(a)\mathcal{A}u\big)d\tau.
\end{align*}
We denote
\begin{align*}
&U_0=-\Lambda^{-2}\nabla\dv\big(e^{\bar\nu\Delta t}u_0\big)+\Lambda^{-2}{\curl}{\curl}\big(e^{\bar\mu\Delta t}u_0\big),\\
&U_1=-\Lambda^{-1}\na h_1+\Lambda^{-1}{\curl}\Omega_1,\\
&U_2=-\Lambda^{-1}\na h_2+\Lambda^{-1}{\curl}\Omega_2,
\end{align*}
where
\begin{equation}\nonumber
\begin{split}
&h_1=-\int_0^te^{{\bar\nu(t-\tau)\Delta }}\Lambda^{-1}\dv(U_0\cdot\na U_0)(\tau)d\tau,\\
&\Omega_1=-\int_0^te^{{\bar\mu(t-\tau)\Delta }}\Lambda^{-1}{\rm{curl}}(U_0\cdot\na U_0)(\tau)d\tau,\\
&h_2=-\int_0^te^{{\bar\nu(t-\tau)\Delta }}\Lambda^{-1}\dv (F_{1}+F_{2}+F_{3})(\tau)d\tau,\\
&\Omega_2=-\int_0^te^{{\bar\mu(t-\tau)\Delta }}\Lambda^{-1}{\rm{curl}} (F_{1}+F_{2})(\tau)d\tau,
\end{split}
\end{equation}
with $F_1, F_2, F_3$ given by
\begin{align*}
&F_{1}=U_0\cdot\na (U_1+U_2)+(U_1+U_2)\cdot\na U_0+(U_1+U_2)\cdot\na(U_1+U_2),
\\&F_{2}=L(a)\mathcal{A}u,
\quad F_{3}=K(a)\nabla a.
\end{align*}
Now we decompose the velocity $u$ as
\beno
u=-\Lambda^{-1}\na h+\Lambda^{-1}{\curl}\Omega=U_0+U_1+U_2.
\eeno

\subsection{The choice of initial data}
Choose $\phi$ as a smooth, radial and non-negative function in $\R^3$ such that
\begin{align}\label{equ:phi}\phi(\xi)=
\left\{
\begin{array}{ll}
1\quad \hbox{for}\quad|\xi|\le 1,\\
0\quad \hbox{for}\quad|\xi|\ge 2.
\end{array}
\right.
\end{align}
Let $N\in\N$ be determined later. We construct the initial data
$(\rho_0, u_0)$  as follows
\begin{align*}
&a_0=\frac{\rho_0}{\bar\rho}-1=\frac{1}{C(N)}\mathcal{F}^{-1}\phi(x),\non\\
&\widehat{u}_0(\xi)=\frac{1}{C(N)}\sum_{k=10}^N 2^{k(1-\frac3p)}\big(\phi(\xi-2^k\tilde{e})+\phi(\xi+2^k\tilde{e}),\,\, i\phi(\xi-2^k\tilde{e})-i\phi(\xi+2^k\tilde{e}),
\,\,0\big),
\end{align*}
where $\tilde{e}=(1,\,1,\,0)$ and $C(N)=2^{\frac
N2(\frac3q-\frac3p+\epsilon)}$ for some $\epsilon>0, q>3, p>6$.
Obviously, the initial velocity $u_0$ is a real-valued function.

\vskip.2cm

The following lemma can be easily verified.

\bthm{Lemma}\label{lem:index}
Let $p>6$. There exist $\epsilon>0$ and $(\widetilde{p},q)$ satisfying
\begin{align*}
&3<q<6,\quad 6<\widetilde{p}<p, \quad\frac3{\widetilde{p}}+\frac3q-1>0,\\
&\max\Big(\frac{2}{\widetilde{p}}-\frac1q-\frac1p,\,\,\frac3{5q}-\frac3{5p}\Big)<\epsilon<\frac13-\frac1q-\frac1p.
\end{align*}
\ethm

Throughout this section, we will fix such a triplet $(\epsilon,
\widetilde{p}, q)$. It is easy to verify that
\begin{align}
&\|a_0\|_{\dot B^{\frac3q}_{q,1}}\le \frac{C}{C(N)},\quad \|u_0\|_{\dot B^{-1+\frac3p}_{p,1}}\le \frac {CN}{C(N)}.\label{eq:u0-p}
\end{align}
Here and in what follows, we denote by $C$  a constant independent of $N$.
Moreover, we have
\begin{equation}
\|u_0\|_{\dot B^{\gamma}_{r,1}}\le
\frac{C2^{N(\gamma-\frac3p+1)}}{C(N)}\non
\end{equation}
for any $r\in[1,\infty],\,\gamma>\f3p-1$. Hence, we get by
Proposition \ref{Prop:heatflow} that
\begin{align}\label{equ:U0}
\|U_0\|_{\widetilde{L}^\rho_T\dot B^{\gamma}_{r,1}}&\le CT^{\frac1{\rho_1}}\|U_0\|_{\widetilde{L}^{\rho_2}_T\dot B^{\gamma}_{r,1}}\le CT^{\frac1{\rho_1}}\|u_0\|_{\dot B^{\gamma-\frac2{\rho_2}}_{r,1}}\non\\
&\le
CT^{\frac1{\rho_1}}\frac{2^{N(\gamma-\frac3p+1-\frac2{\rho_2})}}{C(N)}
\end{align}
for any $r,\rho, \rho_1, \rho_2\in [1,\infty]$ and
$\gamma>\frac3p-1+\frac2{\rho_2}$ with
$\frac1\rho=\frac1{\rho_1}+\frac1{\rho_2}$ .

\subsection{The lower bound estimate of $\|U_1\|_{\dot B^{\frac3p-1}_{p,1}}$}

Since $\dot B^{\frac3p-1}_{p,1}\hookrightarrow \dot
B^{-1}_{\infty,\infty}$, we have
\begin{align}\label{equ:U1_p}
\|U_1(t)\|_{\dot B^{\frac3p-1}_{p,1}}&\ge\|U_1(t)\|_{\dot
B^{-1}_{\infty,\infty}}=\sup_{j\in \Z}2^{-j}\|\Delta_j
U_1(t)\|_{\infty}\ge c\|\Delta_{-4}
U_1(t)\|_{L^\infty}\nonumber\\&=c
\Big\|\int_{\R^3}e^{ix\xi}\varphi(2^4\xi)\widehat{U_1}(t,\xi)d\xi\Big\|_{L^\infty_x}\ge
c\Big|\int_{\R^3}\varphi(2^4\xi) \widehat{U_1}(t,\xi)d\xi\Big|
\end{align}
for some $c>0$ independent of $N$, where $\varphi$ comes from the
Littlewood-Paley decomposition. Set
\begin{align*}
&\mathfrak{U_{11}}=\int_{}\!\int_0^t\varphi(2^4\xi) e^{-{\bar\nu(t-\tau)|\xi|^2}}\mathcal{F}{(U_0\cdot\na U_0)}(\tau, \xi)d\tau d\xi,\\
&\mathfrak{U_{12}}=\int_{}\!\int_0^t\frac{\varphi(2^4\xi)}{|\xi|^{2}}
(e^{-{\bar\mu(t-\tau)|\xi|^2}}-e^{-{\bar\nu(t-\tau)|\xi|^2}})
\mathcal{F}\big({\curl}{\curl}(U_0\cdot\na U_0)\big)(\tau,\xi) d\tau
d\xi,
\end{align*}
we have
\begin{align}\label{equ:U1>}
&\Big|\int_{}\varphi(2^4\xi) \widehat{U_1}(t,\xi)d\xi\Big|\ge
|\mathfrak{U_{11}}+\mathfrak{U_{12}}|.
\end{align}
Due to $\dv\,{\curl} u=0$,  we rewrite $U_0\cdot\na U_0$ as
\begin{align}
U_0\cdot\na U_0=&\frac{1}{2}\nabla\big|\Lambda^{-2}\nabla\dv e^{\bar\nu\Delta \tau}u_0\big|^2+
\dv\big(\Lambda^{-2}{\curl}{\curl}\,e^{\bar\mu\Delta \tau}u_0\otimes U_0\big)\non\\
&-\Lambda^{-2}\na\dv e^{\bar\nu\Delta
\tau}u_0\cdot\na\Lambda^{-2}{\curl}{\curl}\,e^{\bar\mu\Delta
\tau}u_0.\label{eq:U0-d}
\end{align}
This helps us to decompose $\mathfrak{U_{11}}$ into
\begin{align*}
\mathfrak{U_{11}}=\mathfrak{U}_{11}^1+\mathfrak{U}_{11}^2,
\end{align*}
where
\begin{align*}
\mathfrak{U}_{11}^1&=\int\!\!\int_0^t\varphi(2^4\xi) e^{-{\bar\nu(t-\tau)|\xi|^2}}
\mathcal{F}\big(\Lambda^{-2}\na\dv e^{\bar\nu\Delta \tau}u_0\cdot\na\Lambda^{-2}{\curl}{\curl}\,e^{\bar\mu\Delta \tau}u_0\big)(\xi)d\tau d\xi,\\
\mathfrak{U}_{11}^2&=\int\!\!\int_0^t\varphi(2^4\xi)
e^{-{\bar\nu(t-\tau)|\xi|^2}}\Big\{\frac{1}{2}\mathcal{F}
\big(\nabla\big|\Lambda^{-2}\nabla\dv e^{\bar\nu\Delta \tau}u_0\big|^2\big)(\xi)\\
&\qquad\qquad+\mathcal{F}\big(\dv(\Lambda^{-2}{\curl}{\curl}\,e^{\bar\mu\Delta
\tau}u_0\otimes U_0)\big)(\xi)\Big\}d\tau d\xi.
\end{align*}
It follows from \eqref{equ:U1_p} and (\ref{equ:U1>}), we get
\begin{align}\label{equ:U1}
\|U_1(t)\|_{\dot B^{\frac3p-1}_{p,1}}\ge |\mathfrak{U}_{11}^1|-|\mathfrak{U}_{11}^2|-|\mathfrak{U_{12}}|.
\end{align}

In what follows, we consider the case of $t\le 2^{-2N}$.\vspace{0.1cm}

$\bullet$\,\, {\bf The estimate $\mathfrak{U}_{11}^2$}\vspace{0.1cm}

It follows from Lemma \ref{Lem:binesti}, Proposition \ref{Prop:heatflow} and (\ref{eq:u0-p}) that
\begin{align}
|\mathfrak{U}_{11}^2|\le& \Big\|\int_0^te^{-{\bar\nu(t-\tau)\Delta}}\Big\{\f12\nabla\big|\Lambda^{-2}\nabla\dv e^{\bar\nu\Delta \tau}u_0\big|^2+
\dv\big(\Lambda^{-2}{\curl}{\curl}\,e^{\bar\mu\Delta \tau}u_0\otimes U_0\big)\Big\}d\tau\Big\|_{\dot B^{-1}_{\infty,\infty}}\non\\
\le& C\big\|\big|\Lambda^{-2}\nabla\dv e^{\bar\nu\Delta \tau}u_0\big|^2+\Lambda^{-2}{\curl}{\curl}\,e^{\bar\mu\Delta \tau}u_0\otimes U_0\big\|_{L_t^1\dot B^{0}_{\infty,1}}\non\\
\le & C\big\|\big|\Lambda^{-2}\nabla\dv e^{\bar\nu\Delta \tau}u_0\big|^2+\Lambda^{-2}{\curl}{\curl}\,e^{\bar\mu\Delta \tau}u_0\otimes U_0\big\|_{L_t^1\dot B^{\f3p}_{p,1}}\non\\
\le &C\|e^{\bar\nu\Delta \tau}u_0\|_{\widetilde{L}_t^2\dot B^{\frac3p}_{p,1}}^2+C\|e^{\bar\mu\Delta \tau}u_0\|_{\widetilde{L}_t^2\dot B^{\frac3p}_{p,1}}^2\non\\
\le &C\|u_0\|_{\dot B^{\frac3p-1}_{p,1}}^2\le \frac {CN^2}{C(N)^2}.\label{eq:U11-2}
\end{align}

$\bullet$\,\, {\bf The estimate of $\mathfrak{U_{12}}$}\vspace{0.1cm}

Using (\ref{eq:U0-d}), we decompose $\mathfrak{U_{12}}$ as
\begin{align}
\mathfrak{U}_{12}=\mathfrak{U}_{12}^1+\mathfrak{U}_{12}^{2},\non
\end{align}
where
\begin{align*}
\mathfrak{U}_{12}^{1}=&\int_{}\!\!\int_0^t\varphi(2^4\xi)
(e^{-{\bar\mu(t-\tau)|\xi|^2}}-e^{-{\bar\nu(t-\tau)|\xi|^2}}){|\xi|^{-2}}
\nonumber\\&\quad\times\mathcal{F}\big({{\curl}{\curl}(\Lambda^{-2}\na\dv
e^{\bar\nu\Delta \tau}u_0
\cdot\na\Lambda^{-2}{\curl}{\curl}\,e^{\bar\mu\Delta \tau}u_0)\big)}(\xi)d\tau d\xi,\\
\mathfrak{U}_{12}^{2}=&\int_{}\!\!\int_0^t\varphi(2^4\xi) (e^{-{\bar\mu(t-\tau)|\xi|^2}}-e^{-{\bar\nu(t-\tau)|\xi|^2}}){|\xi|^{-2}}
\mathcal{F}\big({\curl}{\curl}\big\{\frac{1}{2}\nabla\big|\Lambda^{-2}\nabla\dv e^{\bar\nu\Delta \tau}u_0\big|^2\nonumber\\
&\qquad+\dv(\Lambda^{-2}{\curl}{\curl}\,e^{\bar\mu\Delta
\tau}u_0\otimes U_0)\big\}\big)(\xi)d\tau d\xi.
\end{align*}

By the same argument as the one deriving to (\ref{eq:U11-2}), we
infer
\begin{align}\label{eq:U12-2}
|\mathfrak{U}_{12}^{2}|\le C\|u_0\|_{\dot B^{\frac3p-1}_{p,1}}^2\le \frac {CN^2}{C(N)^2}.
\end{align}
We denote by $a(\xi)$ the symbol of the operator $\curl\curl$. Then $\mathfrak{U}_{12}^{1}$ is written as
\begin{align*}
&\int_{}\!\!\int_0^t\varphi(2^4\xi) {|\xi|^{-2}}
(e^{-{\bar\mu(t-\tau)|\xi|^2}}-e^{-{\bar\nu(t-\tau)|\xi|^2}})\int
\frac{e^{-\bar\nu\tau|\eta|^2}}{|\eta|^2}\frac{e^{-\bar\mu\tau|\xi-\eta|^2}}{|\xi-\eta|^2}
\nonumber\\& \qquad\times{a(\xi)}\big(\mathcal{F}(\na\dv
u_0)(\eta)\cdot\mathcal{F}(\na{\curl}{\curl}u_0)(\xi-\eta)\big)d\eta
d\tau d\xi.
\end{align*}
Due to the choice of $u_0$, we find that $|\eta|\gg |\xi|\sim 1$,
this yields that \beno
&&\int_0^t (e^{-{\bar\mu(t-\tau)|\xi|^2}}-e^{-{\bar\nu(t-\tau)|\xi|^2}})e^{-\bar\nu\tau|\eta|^2-\bar\mu\tau|\xi-\eta|^2}d\tau\\
&&=\frac{e^{-\bar\mu t|\xi|^2}-e^{-\bar\nu t|\eta|^2-\bar\mu t|\xi-\eta|^2}}{\bar\nu|\eta|^2+\bar\mu|\xi-\eta|^2-\bar\mu|\xi|^2}-
\frac{e^{-\bar\nu t|\xi|^2}-e^{-\bar\nu t|\eta|^2-\bar\mu t|\xi-\eta|^2}}{\bar\nu|\eta|^2+\bar\mu|\xi-\eta|^2-\bar\nu|\xi|^2}\\
&&\le Ct^2|\eta|^{2}.
\eeno
Then we obtain
\begin{align}
|\mathfrak{U}_{12}^{1}|&\le \frac{C t^2}{C(N)^2}\sum_{k=10}^N
2^{(5-\f6p)k}\le \frac{C 2^{(5-\f6p)N}t^2}{C(N)^2}.\non
\end{align}
This along with (\ref{eq:U12-2}) implies  that
\begin{align}\label{eq:U12}
|\mathfrak{U}_{12}|\le \frac
{CN^2}{C(N)^2}+\frac{C2^{(5-\f6p)N}t^2}{C(N)^2}.
\end{align}

$\bullet$\,\, {\bf The estimate of $\mathfrak{U}_{11}^1$}\vspace{0.1cm}

The $\ell'$-th component $\mathfrak{U}_{11}^{1\ell'} $ of $\mathfrak{U}_{11}^1$ is given by
\begin{align}
\mathfrak{U}_{11}^{1,\ell'} =&\int\!\!\int_0^t\varphi(2^4\xi) e^{-{\bar\nu(t-\tau)|\xi|^2}}i\int
\frac{\eta_\ell\eta_{m}}{|\eta|^2|\xi-\eta|^2}e^{-{\bar\nu\tau|\eta|^2}}\widehat{u}_{0{m}}(\eta)
(\xi-\eta)_{\ell}e^{-{\bar\mu\tau|\xi-\eta|^2}}(\xi-\eta)_{m'}
\nonumber\\&\qquad\times\big((\xi-\eta)_{m'}\widehat{u}_{0{\ell'}}(\xi-\eta)-(\xi-\eta)_{\ell'}\widehat{u}_{0{m'}}(\xi-\eta)\big)
d\eta d\tau d\xi\non\\
\triangleq &\int\!\!\int_0^t\varphi(2^4\xi)
e^{-{\bar\nu(t-\tau)|\xi|^2}}iA(t,\xi)d\tau d\xi.\non
\end{align}
In what follows, we consider the case of $\ell'=1$. From the
construction of $u_0$, we find that \beno
A(\tau,\xi)=\frac{1}{C(N)^2}\sum_{k=10}^N
2^{2k(1-\frac3p)}\big(A_1+A_2), \eeno where
\begin{align*}
A_1=\int\!\big[&\eta_1(\xi_2-\eta_2)^2+\eta_1(\xi_3-\eta_3)^2-\eta_2(\xi_1-\eta_1)(\xi_2-\eta_2)
+i\eta_1(\xi_1-\eta_1)(\xi_2-\eta_2)\\&+i\eta_2(\xi_2-\eta_2)^2+i\eta_2(\xi_3-\eta_3)^2\big]
\frac{\eta_\ell(\xi_{\ell}-\eta_{\ell})}{|\eta|^2|\xi-\eta|^2}
e^{-{\bar\nu\tau|\eta|^2}-{\bar\mu\tau|\xi-\eta|^2}}
\phi(\eta-2^k\tilde{e})\phi(\xi-\eta+2^k\tilde{e})d\eta,\\
A_2=\int\!\big[&\eta_1(\xi_2-\eta_2)^2+\eta_1(\xi_3-\eta_3)^2-\eta_2(\xi_1-\eta_1)(\xi_2-\eta_2)-i\eta_1(\xi_1-\eta_1)(\xi_2-\eta_2)
\nonumber\\&- i\eta_2(\xi_2-\eta_2)^2-i\eta_2(\xi_3-\eta_3)^2\big]
\frac{\eta_\ell(\xi_{\ell}-\eta_{\ell})}{|\eta|^2|\xi-\eta|^2}
e^{-{\bar\nu\tau|\eta|^2}-{\bar\mu\tau|\xi-\eta|^2}}
\phi(\eta+2^k\tilde{e})\phi(\xi-\eta-2^k\tilde{e})d\eta.
\end{align*}
Making a change of variable, we obtain
\begin{align*}
A_1=\int\big[&(\eta_1+2^k)((\xi_2-\eta_2-2^k)^2+(\xi_3-\eta_3)^2)-(\eta_2+2^k)(\xi_1-\eta_1-2^k)(\xi_2-\eta_2-2^k)
\nonumber\\&+i(\eta_1+2^k)(\xi_1-\eta_1-2^k)(\xi_2-\eta_2-2^k)+i(\eta_2+2^k)((\xi_2-\eta_2-2^k)^2+(\xi_3-\eta_3)^2)\big]
\nonumber\\&\times\frac{(\eta+2^k\tilde{e})_\ell(\xi-\eta-2^k\tilde{e})_{\ell}}
{|\eta+2^k\tilde{e}|^2|\xi-\eta-2^k\tilde{e}|^2}
e^{-{\bar\nu\tau|\eta+2^k\tilde{e}|^2}-{\bar\mu\tau|\xi-\eta-2^k\tilde{e}|^2}}\phi(\eta)\phi(\xi-\eta)d\eta,
\end{align*}
and
\begin{align*}
A_2=\int\big[&(\eta_1-2^k)((\xi_2-\eta_2+2^k)^2+(\xi_3-\eta_3)^2)-(\eta_2-2^k)(\xi_1-\eta_1+2^k)(\xi_2-\eta_2+2^k)
\nonumber\\&-i(\eta_1-2^k)(\xi_1-\eta_1+2^k)(\xi_2-\eta_2+2^k)-i(\eta_2-2^k)((\xi_2-\eta_2+2^k)^2+(\xi_3-\eta_3)^2)
\big]\nonumber\\&\times\frac{(\eta-2^k\tilde{e})_\ell(\xi-\eta+2^k\tilde{e})_{\ell}}
{|\eta-2^k\tilde{e}|^2|\xi-\eta+2^k\tilde{e}|^2}
e^{-{\bar\nu\tau|\eta-2^k\tilde{e}|^2}-{\bar\mu\tau|\xi-\eta+2^k\tilde{e}|^2}}\phi(\eta)\phi(\xi-\eta)d\eta.
\end{align*}
Due to the choice of $\phi$, we find
\begin{align*}
A_1&=\int
\big(-i2^{2+5k}+O(2^{4k})\big)\frac{e^{-{\bar\nu\tau|\eta+2^k\tilde{e}|^2}-{\bar\mu\tau|\xi-\eta-2^k\tilde{e}|^2}}}
{|\eta+2^k\tilde{e}|^2|\xi-\eta-2^k\tilde{e}|^2}
\phi(\eta)\phi(\xi-\eta)d\eta,\\
A_2&=\int\big(-i2^{2+5k}+O(2^{4k})\big)\frac{e^{-{\bar\nu\tau|\eta-2^k\tilde{e}|^2}-{\bar\mu\tau|\xi-\eta+2^k\tilde{e}|^2}}}
{|\eta-2^k\tilde{e}|^2|\xi-\eta+2^k\tilde{e}|^2}
\phi(\eta)\phi(\xi-\eta)d\eta.
\end{align*}
This yields that
\begin{align}
\mathfrak{U}_{11}^{1,1}=\mathfrak{U}_{11}^{11,1}+\mathfrak{U}_{11}^{12,1},\non
\end{align}
where
\begin{align*}
\mathfrak{U}_{11}^{11,1}=\frac{1}{C(N)^2}&\sum_{k=10}^N
2^{2k(1-\frac3p)}\int\!\!\!\int_0^t\varphi(2^4\xi)
e^{-{\bar\nu(t-\tau)|\xi|^2}} \int
2^{2+5k}\phi(\eta)\phi(\xi-\eta)\nonumber\\&\times
\bigg\{\frac{e^{-{\bar\nu\tau|\eta+2^k\tilde{e}|^2}-{\bar\mu\tau|\xi-\eta-2^k\tilde{e}|^2}}}
{|\eta+2^k\tilde{e}|^2|\xi-\eta-2^k\tilde{e}|^2}
+\frac{e^{-{\bar\nu\tau|\eta-2^k\tilde{e}|^2}-{\bar\mu\tau|\xi-\eta+2^k\tilde{e}|^2}}}
{|\eta-2^k\tilde{e}|^2|\xi-\eta+2^k\tilde{e}|^2}\bigg\}d\eta d\tau d\xi,\\
\mathfrak{U}_{11}^{12,1}=\frac{1}{C(N)^2}&\sum_{k=10}^N
2^{2k(1-\frac3p)}\int\!\!\!\int_0^t\varphi(2^4\xi)
e^{-{\bar\nu(t-\tau)|\xi|^2}} i\int
O(2^{4k})\phi(\eta)\phi(\xi-\eta)\nonumber\\&\times
\bigg\{\frac{e^{-{\bar\nu\tau|\eta+2^k\tilde{e}|^2}-{\bar\mu\tau|\xi-\eta-2^k\tilde{e}|^2}}}
{|\eta+2^k\tilde{e}|^2|\xi-\eta-2^k\tilde{e}|^2}
+\frac{e^{-{\bar\nu\tau|\eta-2^k\tilde{e}|^2}-{\bar\mu\tau|\xi-\eta+2^k\tilde{e}|^2}}}
{|\eta-2^k\tilde{e}|^2|\xi-\eta+2^k\tilde{e}|^2}\bigg\}d\eta d\tau
d\xi.
\end{align*}
After integrating with respect to $\tau$, we get
\begin{align*}
\mathfrak{U}_{11}^{11,1}=&\frac{4}{C(N)^2}\sum_{k=10}^N
2^{2k(1-\frac3p)}\int\!\!\varphi(2^4\xi) \int
2^{5k}\phi(\eta)\phi(\xi-\eta)\nonumber\\&\times
\bigg\{\frac{e^{-{\bar\nu t|\xi|^2}}-e^{-{\bar\nu
t|\eta+2^k\tilde{e}|^2}-{\bar\mu t|\xi-\eta-2^k\tilde{e}|^2}}}
{|\eta+2^k\tilde{e}|^2|\xi-\eta-2^k\tilde{e}|^2(\bar\nu
|\eta+2^k\tilde{e}|^2+\bar\mu
|\xi-\eta-2^k\tilde{e}|^2-\bar\nu|\xi|^2)}\nonumber\\&
\qquad+\frac{e^{-{\bar\nu t|\xi|^2}}-e^{-{\bar\nu
t|\eta-2^k\tilde{e}|^2}-{\bar\mu t|\xi-\eta+2^k\tilde{e}|^2}}}
{|\eta-2^k\tilde{e}|^2|\xi-\eta+2^k\tilde{e}|^2(\bar\nu
|\eta-2^k\tilde{e}|^2+\bar\mu
|\xi-\eta+2^k\tilde{e}|^2-\bar\nu|\xi|^2)}\bigg\}d\eta d\xi.
\end{align*}
Using Taylor's formula, we infer
\beno
&&\frac{e^{-{\bar\nu t|\xi|^2}}-e^{-{\bar\nu t|\eta+2^k\tilde{e}|^2}-{\bar\mu t|\xi-\eta-2^k\tilde{e}|^2}}}
{\bar\nu |\eta+2^k\tilde{e}|^2+\bar\mu |\xi-\eta-2^k\tilde{e}|^2-\bar\nu|\xi|^2}=t
+O(t^22^{2k}),\\
&&\frac{e^{-{\bar\nu t|\xi|^2}}-e^{-{\bar\nu t|\eta-2^k\tilde{e}|^2}-{\bar\mu t|\xi-\eta+2^k\tilde{e}|^2}}}
{\bar\nu |\eta-2^k\tilde{e}|^2+\bar\mu |\xi-\eta+2^k\tilde{e}|^2-\bar\nu|\xi|^2}=t
+O(t^22^{2k}),
\eeno
from which, it follows that
\begin{align*}
\mathfrak{U}_{11}^{11,1}=&\frac{4}{C(N)^2}\sum_{k=10}^N
2^{2k(1-\frac3p)}\int\!\!\int\varphi(2^4\xi)
2^{5k}\phi(\eta)\phi(\xi-\eta)\nonumber\\&\times\bigg\{\frac{t}
{|\eta+2^k\tilde{e}|^2|\xi-\eta-2^k\tilde{e}|^2}
+\frac{t}{|\eta-2^k\tilde{e}|^2|\xi-\eta+2^k\tilde{e}|^2}+O(t^22^{-2k})\bigg\}d\eta
d\xi,
\end{align*}
hence,
\begin{align}
\mathfrak{U}_{11}^{11,1}
\ge&\frac{ct}{C(N)^2}\sum_{k=10}^N 2^{3k-\frac6pk}-\frac{Ct^2}{C(N)^2}\sum_{k=10}^N2^{5k-\f6pk}\non\\
\ge&\frac{ct2^{3N-\frac6pN}}{C(N)^2}-\frac{Ct^22^{5N-\f6pN}}{C(N)^2}.\label{eq:U-11-111}
\end{align}
Similarly, we can deduce
\begin{align}\label{eq:U-11-112}
|\mathfrak{U}_{11}^{12,1}|\le &\frac{C}{C(N)^2}\sum_{k=10}^N
\int\!\!\int2^{2k(1-\frac3p)}\bigg\{\frac{e^{-{\bar\nu
t|\xi|^2}}-e^{-{\bar\nu t|\eta+2^k\tilde{e}|^2}-{\bar\mu
t|\xi-\eta-2^k\tilde{e}|^2}}}{\bar\nu |\eta+2^k\tilde{e}|^2+\bar\mu
|\xi-\eta-2^k\tilde{e}|^2-\bar\nu|\xi|^2}\nonumber\\&\quad+\frac{e^{-{\bar\nu
t|\xi|^2}}-e^{-{\bar\nu t|\eta-2^k\tilde{e}|^2}-{\bar\mu
t|\xi-\eta+2^k\tilde{e}|^2}}}{\bar\nu |\eta-2^k\tilde{e}|^2+\bar\mu
|\xi-\eta+2^k\tilde{e}|^2-\bar\nu|\xi|^2}\bigg\}
\phi(\eta)\phi(\xi-\eta)d\eta
d\xi\nonumber\\\le&\frac{C}{C(N)^2}\sum_{k=10}^N 2^{-\frac6pk}\le
\frac {C}{C(N)^2}.
\end{align}
Hence, we conclude that \ben\label{eq:U-11-1}
|\mathfrak{U}_{11}^{1}|\ge |\mathfrak{U}_{11}^{11}|\ge
\frac{ct2^{3N-\frac6pN}}{C(N)^2}-\frac{Ct^22^{5N-\f6pN}}{C(N)^2}-\frac
{C}{C(N)^2}. \een

Summing up (\ref{equ:U1}), (\ref{eq:U11-2}), (\ref{eq:U12}) and (\ref{eq:U-11-1}), we obtain
\begin{align}
\|U_1(t)\|_{\dot B^{\frac3p-1}_{p,1}}\ge \frac{ct2^{3N-\frac6pN}}{C(N)^2}-\frac{Ct^22^{5N-\f6pN}}{C(N)^2}-\frac {CN^2}{C(N)^2}.\non
\end{align}
Choosing $t=2^{-2(1+\epsilon)N}$ and recalling $C(N)^2=2^{N(\f3q-\f3p+\varepsilon)}$, we get
\ben\label{eq:U1-final}
\|U_1(t)\|_{\dot B^{\frac3p-1}_{p,1}}\ge  c2^{(1-\f3q-\f3p-3\epsilon)N}
\een
for some $c>0$ independent of $N$.

\subsection{The estimate of $\|U_1\|_{L^1_T\dot B^{\frac3q+1}_{q,1}\cap \widetilde{L}^2_T\dot B^{\frac3q}_{q,1}}$}

Let $(p,q,\widetilde{p})$ be given as in Lemma \ref{lem:index}. It
follows from H\"{o}lder's inequality and Proposition
\ref{Prop:heatflow} that
\begin{align}\label{equ:U1estimate0}
&\|U_1\|_{{L}^1_T\dot B^{\frac3q+1}_{q,1}}+\|U_1\|_{\widetilde{L}^2_T\dot B^{\frac3q}_{q,1}}\non\\
&\le C\big(\|h_1\|_{L^1_T\dot B^{\frac3q+1}_{q,1}}+\|h_1\|_{\widetilde{L}^2_T\dot B^{\frac3q}_{q,1}}+\|\Omega_1\|_{{L}^1_T\dot B^{\frac3q+1}_{q,1}}+\|\Omega_1\|_{\widetilde{L}^2_T\dot B^{\frac3q}_{q,1}}\big)\nonumber\\&
\le CT^{\frac12}\Big(\|h_1\|_{\widetilde{L}^2_T\dot B^{\frac3q+1}_{q,1}}+\|\Omega_1\|_{\widetilde{L}^2_T\dot B^{\frac3q+1}_{q,1}}
+\|h_1\|_{\widetilde{L}^\infty_T\dot B^{\frac3q}_{q,1}}+\|\Omega_1\|_{\widetilde{L}^\infty_T\dot B^{\frac3q}_{q,1}}\Big)\nonumber\\&\le CT^{\frac12}\|U_0\cdot\na U_0\|_{L^1_T\dot B^{\frac3q}_{q,1}}.
\end{align}
We infer from Lemma \ref{Lem:binestis>0}, $\dot
B^{\frac3p}_{p,1}\hookrightarrow L^\infty $ and \eqref{equ:U0} that
\begin{align*}
\|U_0\cdot\na U_0\|_{L^1_T\dot B^{\frac3q}_{q,1}}&\le C\big(\|U_0\|_{L^\infty_TL^\infty}\|\na U_0\|_{L^1_T\dot B^{\frac3q}_{q,1}}+\|U_0\|_{L^1_T\dot B^{\frac3q}_{q,1}}\|\na U_0\|_{L^\infty_TL^\infty}\big)
\nonumber\\&\le CT^\frac12\big(\|U_0\|_{L^\infty_T\dot B^{\frac3p}_{p,1}}\|U_0\|_{\widetilde{L}^2_T\dot B^{\frac3q+1}_{q,1}}+\|U_0\|_{\widetilde{L}^2_T\dot B^{\frac3q}_{q,1}}\|U_0\|_{L^\infty_T\dot B^{\frac3p+1}_{p,1}}\big)\nonumber\\
&\le CT^\frac12\frac{2^{N(\frac3q-\frac3p+2)}}{C(N)^2},
\end{align*}
from which and \eqref{equ:U1estimate0}, it follows that
\begin{align}\label{equ:U1estimateq}
\|U_1\|_{{L}^1_T\dot B^{\frac3q+1}_{q,1}}+\|U_1\|_{\widetilde{L}^2_T\dot B^{\frac3q}_{q,1}}\le CT\frac {2^{N(\frac3q-\frac3p+2)}}{C(N)^2}.
\end{align}
Similarly, we have
\begin{align}
&\|U_1\|_{L^1_T\dot B^{\frac3p+1}_{p,1}}+\|U_1\|_{\widetilde{L}^2_T\dot B^{\frac3p}_{p,1}}\le CT\frac {2^{2N}}{C(N)^2},\label{equ:U1estimatep}\\
&\|U_1\|_{L^1_T\dot B^{\frac3{\widetilde{p}}+1}_{\widetilde{p},1}}
+\|U_1\|_{\widetilde{L}^2_T\dot B^{\frac3{\widetilde{p}}}_{\widetilde{p},1}}\le CT\frac {2^{N(\frac3{\widetilde{p}}-\frac3p+2)}}{C(N)^2}.\label{equ:U1estimatetildep}
\end{align}

\subsection{The estimates of $F_i(i=1,2,3)$}
Recalling  that
\begin{align*}
F_3=K(a)\nabla a=(K(a)-1)\nabla a+\nabla a.
\end{align*}
Then it follows from Lemma \ref{Lem:binesti} and Lemma
\ref{Lem:nonesti} that
\begin{align}\label{equ:F3}
\|F_3\|_{L^1_T\dot B^{\frac3q-1}_{q,1}}&\le CT\|K(a)-1\|_{L^\infty_T\dot B^{\frac3q}_{q,1}}\|\nabla {a}\|_{L^\infty_T\dot B^{\frac3q-1}_{q,1}}+\|a\|_{L^1_T\dot B^{\frac3q}_{q,1}}\nonumber\\
&\le
CT\big(1+\|a\|_{L^\infty_T(L^\infty)}\big)^2\|a\|_{L^\infty_T\dot
B^{\frac3q}_{q,1}}^2+CT\|a\|_{L^\infty_T\dot B^{\frac3q}_{q,1}}.
\end{align}
We write $F_2$ as
\begin{align*}
F_2=L(a)\mathcal{A}u =L(a)\mathcal{A}(U_0+U_1)+L(a)\mathcal{A}U_2.
\end{align*}
By making use of Lemma \ref{Lem:binesti} again, (\ref{equ:U0}) and
(\ref{equ:U1estimateq}), we obtain
\begin{align}
\|L(a)\mathcal{A} (U_0+U_1)\|_{L^1_T\dot B^{\frac3q-1}_{q,1}}
&\le C\|L(a)\|_{L^\infty_T\dot B^{\frac3q}_{q,1}}\|\mathcal{A} (U_0+U_1)\|_{L^1_T\dot B^{\frac3q-1}_{q,1}}
\nonumber\\&\le C\big(\|a\|_{L^\infty_T(L^\infty)}+1\big)^2\|a\|_{L^\infty_T\dot B^{\frac3q}_{q,1}}\|U_0+U_1\|_{L^1_T\dot B^{\frac3q+1}_{q,1}}
\nonumber\\&\le C\big(\|a\|_{L^\infty_T(L^\infty)}+1\big)^2\|a\|_{L^\infty_T\dot B^{\frac3q}_{q,1}}T\frac{2^{N(\frac{3}{q}-\frac3p+2)}}{C(N)},\non
\end{align}
and
\begin{align}\label{equ:Fy2part3}
\|L(a)\mathcal{A}U_2\|_{L^1_T\dot B^{\frac3q-1}_{q,1}}&\le
C\|L(a)\|_{L^\infty_T\dot B^{\frac3q}_{q,1}}\|\mathcal{A}U_2\|_{L^1_T\dot B^{\frac3q-1}_{q,1}}
\nonumber\\&\le C\big(\|a\|_{L^\infty_T(L^\infty)}+1\big)^2\|a\|_{L^\infty_T\dot B^{\frac3q}_{q,1}}\|U_2\|_{L^1_T\dot B^{\frac3q+1}_{q,1}}.\non
\end{align}
This above  two estimates give that \ben\label{equ:F2}
\|F_2\|_{L^1_T\dot B^{\frac3q-1}_{q,1}}\le
C\big(\|a\|_{L^\infty_T(L^\infty)}+1\big)^2\|a\|_{L^\infty_T\dot
B^{\frac3q}_{q,1}}
\Big(T\frac{2^{N(\frac{3}{q}-\frac3p+2)}}{C(N)}+\|U_2\|_{L^1_T\dot
B^{\frac3q+1}_{q,1}}\Big). \een

Now let us turn to the estimate of $F_1$, which is given by \beno
F_{1}=U_0\cdot\na (U_1+U_2)+(U_1+U_2)\cdot\na
U_0+(U_1+U_2)\cdot\na(U_1+U_2). \eeno Due to $\f 3 p-\f3q+1> 0$ and
$\f 3 {\widetilde{p}}+\f3q-1>0$,  we apply (\ref{Bonydecom}), Lemma
\ref{lem:biest-n}, (\ref{equ:U0}) and (\ref{equ:U1estimateq}) to get
\ben\label{equ:F1-1}
&&\|U_0\cdot\na (U_1+U_2)\|_{L^1_T\dot B^{\frac3q-1}_{q,1}}\nonumber\\
&&\le C\|U_0\|_{\widetilde{L}^2_T\dot B^{\frac3p}_{p,1}}
\|\na(U_1+U_2)\|_{\widetilde{L}^2_T\dot B^{\frac3q-1}_{q,1}}+C\|U_0\|_{\widetilde{L}^2_T\dot B^{\frac3{\widetilde{p}}}_{\widetilde{p},1}}
\|\na(U_1+U_2)\|_{\widetilde{L}^2_T\dot B^{\frac3q-1}_{q,1}}\nonumber\\
&&\le CT^{\frac12}\Big(\frac{2^{N}}{C(N)}+\frac{2^{N(\frac{3}{\widetilde{p}}-\frac3p+1)}}{C(N)}\Big)\Big(T\frac {2^{N(\frac3q-\frac3p+2)}}{C(N)^2}+\|U_2\|_{\widetilde{L}^2_T\dot B^{\frac3q}_{q,1}}\Big),
\een
and
\begin{align}
\|U_1\cdot\na U_0\|_{L^1_T\dot B^{\frac3q-1}_{q,1}}\le& C\Big(\|\na U_0\|_{\widetilde{L}^2_T\dot B^{\frac3p-1}_{p,1}}\|U_1\|_{\widetilde{L}^2_T\dot B^{\frac3q}_{q,1}}+\|\na U_0\|_{\widetilde{L}^2_T\dot B^{\frac3{q}-1}_{q,1}}\|U_1\|_{\widetilde{L}^2_T\dot B^{\frac3p}_{p,1}}\nonumber\\
&\qquad+\|\na U_0\|_{\widetilde{L}^2_T\dot B^{\frac3{\widetilde{p}}-1}_{\widetilde{p},1}}\|U_1\|_{\widetilde{L}^2_T\dot B^{\frac3q}_{q,1}}\Big)\nonumber\\
\le& CT^{\frac32}\Big(\frac {2^{N(\frac3q-\frac3p+3)}}{C(N)^3}+\frac {2^{N(\frac3{\widetilde{p}}+\frac3q-\frac6p+3)}}{C(N)^3}\Big),
\end{align}
and
\begin{align}
\|&U_2\cdot\na U_0\|_{L^1_T\dot B^{\frac3q-1}_{q,1}}\le C\Big(\|\na U_0\|_{\widetilde{L}^2_T\dot B^{\frac3p-1}_{p,1}}\|U_2\|_{\widetilde{L}^2_T\dot B^{\frac3q}_{q,1}}+\|\na U_0\|_{L^1_T\dot B^{\frac3{p}}_{p,1}}\|U_2\|_{\widetilde{L}^\infty_T\dot B^{\frac3q-1}_{q,1}}\nonumber\\
&\qquad\qquad\qquad\qquad\qquad+
\|\na U_0\|_{\widetilde{L}^2_T\dot B^{\frac3{\widetilde{p}}-1}_{\widetilde{p},1}}\|U_2\|_{\widetilde{L}^2_T\dot B^{\frac3q}_{q,1}}\Big)\nonumber\\
&\le CT^{\frac12}\Big(\frac {2^{N}}{C(N)}+\frac
{2^{N(\frac3{\widetilde{p}}-\frac3p+1)}}{C(N)}\Big)\|U_2\|_{\widetilde{L}^2_T\dot
B^{\frac3q}_{q,1}}+CT\frac
{2^{2N}}{C(N)}\|U_2\|_{\widetilde{L}^\infty_T\dot
B^{\frac3q-1}_{q,1}}.
\end{align}
By Lemma \ref{lem:biest-n} and
(\ref{equ:U1estimateq})-(\ref{equ:U1estimatetildep}), we have
\begin{align}
&\|U_1\cdot\na (U_1+U_2)\|_{L^1_T\dot B^{\frac3q-1}_{q,1}}+\|U_2\cdot\na U_1\|_{L^1_T\dot B^{\frac3q-1}_{q,1}}\nonumber\\
&\le C\|U_1\|_{\widetilde{L}^2_T\dot B^{\frac3p}_{p,1}}\|\na (U_1+U_2)\|_{\widetilde{L}^2_T\dot B^{\frac3q-1}_{q,1}}
+C\|U_1\|_{\widetilde{L}^2_T\dot B^{\frac3{\widetilde{p}}}_{\widetilde{p},1}}\|\na (U_1+U_2)\|_{\widetilde{L}^2_T\dot B^{\frac3q-1}_{q,1}}\nonumber\\
&\qquad+C\|\nabla U_1\|_{L^1_T\dot B^{\frac3p}_{p,1}}\| U_2\|_{\widetilde{L}^\infty_T\dot B^{\frac3q-1}_{q,1}}+\|\nabla U_1\|_{L^1_T\dot B^{\frac3{\widetilde{p}}}_{\widetilde{p},1}}\| U_2\|_{\widetilde{L}^\infty_T\dot B^{\frac3q-1}_{q,1}}\non\\
&\le  CT^{2}\Big(\frac {2^{N(\frac3q-\frac3p+4)}}{C(N)^4}+\frac {2^{N(\frac3{\widetilde{p}}+\frac3q-\frac6p+4)}}{C(N)^4}\Big)
+CT\Big(\frac{2^{2N}}{C(N)^2}+\frac{2^{N(\frac3{\widetilde{p}}-\frac3p+2)}}{C(N)^2}\Big)\|U_2\|_{\widetilde{L}^2_T\dot B^{\frac3q}_{q,1}}\non\\
&\qquad+CT\Big(\frac{2^{2N}}{C(N)^2}+\frac{2^{N(\frac3{\widetilde{p}}-\frac3p+2)}}{C(N)^2}\Big)\|U_2\|_{\widetilde{L}^\infty_T\dot
B^{\frac3q-1}_{q,1}}.
\end{align}
We infer from Lemma \ref{Lem:binesti} that
\begin{align}\label{equ:F1-5}
&\|U_2\cdot\na U_2\|_{L^1_T\dot B^{\frac3q-1}_{q,1}}\le C\|U_2\|_{\widetilde{L}^2_T\dot B^{\frac3q}_{q,1}}\|\na U_2\|_{\widetilde{L}^2_T\dot B^{\frac3q-1}_{q,1}}\le C\|U_2\|^2_{\widetilde{L}^2_T\dot B^{\frac3q}_{q,1}}.
\end{align}

Summing up (\ref{equ:F1-1})-(\ref{equ:F1-5}), we obtain
\begin{align}\label{equ:F1}
&\|F_1\|_{L^1_T\dot B^{\frac3q-1}_{q,1}}\le \f {CT^\f32} {C(N)^3}\big(2^{N(\f3q-\f3p+3)}+2^{N(\f 3 {\widetilde{p}}+\f3q-\f6p+3)}\big)\non\\
&\quad+\f {CT^2} {C(N)^4}\big(2^{N(\f3q-\f3p+4)}+2^{N(\f 3 {\widetilde{p}}+\f3q-\f6p+4)}\big)
+C\|U_2\|^2_{\widetilde{L}^2_T\dot B^{\frac3q}_{q,1}}\non\\
&\quad+\Big\{\f {CT^\f12}{C(N)}\big(2^N+2^{N(\f 3 {\widetilde{p}}-\f3p+1)}\big)
+CT\Big(\frac{2^{2N}}{C(N)}+\frac{2^{N(\frac3{\widetilde{p}}-\frac3p+2)}}{C(N)^2}\Big)\Big\}\|U_2\|_{\widetilde{L}^2_T\dot B^{\frac3q}_{q,1}}\non\\
&\quad+\f {CT}
{C(N)}\Big(\frac{2^{2N}}{C(N)}+\frac{2^{N(\frac3{\widetilde{p}}-\frac3p+2)}}{C(N)}+{2^{2N}}\Big)
\|U_2\|_{\widetilde{L}^\infty_T\dot B^{\frac3q-1}_{q,1}}.
\end{align}

\subsection{Proof of Theorem \ref{thm:illposed-bNS}}
We denote
$$
X_T=\|{a}\|_{\widetilde{L}^\infty_T\dot B^{\frac3q}_{q,1}},\quad Y_T=\|U_2\|_{\widetilde{L}^\infty_T\dot B^{\frac3q-1}_{q,1}}+\|U_2\|_{L^1_T\dot B^{\frac3q+1}_{q,1}}.
$$
For $T\le T_0=2^{-2(1+\epsilon)N}$, it follows from Proposition
\ref{Prop:heatflow}, (\ref{equ:F3}), (\ref{equ:F2}) and
(\ref{equ:F1}) that
\begin{align}
Y_T\le& C\sum_{i=1}^3\|F_i\|_{L^1_T\dot B^{\frac3q-1}_{q,1}}\non\\
\le & \f {CT^\f322^{N(\f 3 {\widetilde{p}}+\f3q-\f6p+3)}} {C(N)^3}+\f {CT^\f122^{N(\f 3 {\widetilde{p}}-\f3p+1)}}{C(N)}Y_T+CY_T^2\non\\
&+\f {CT2^{N(\f3q-\f3p+2)}}
{C(N)}(1+X_T)^3X_T+C(1+X_T)^2X_TY_T.\label{equ:YT}
\end{align}
On the other hand, we infer from Proposition \ref{Prop:transport} that
\begin{align}\label{equ:tildeaestimate}
&X_T \le C \exp(\|\na u\|_{L^1_T\dot B^{\frac3p}_{p,1}})
\big(\|{a}_0\|_{\dot B^{\frac3q}_{q,1}}+\|\dv u+a\dv u\|_{L^1_T\dot B^{\frac3q}_{q,1}}\big).\non
\end{align}
Thanks to Lemma \ref{Lem:binesti}, (\ref{equ:U0}) and (\ref{equ:U1estimateq}),  we have
\begin{align}
\|a\dv u\|_{L^1_T\dot B^{\frac3q}_{q,1}}&\le C\|\dv(U_0+U_1)\|_{L^1_T\dot B^{\frac3q}_{q,1}}X_T+CY_TX_T\nonumber\\&\le
\frac{CT2^{N(\frac3{q}-\frac3p+2)}}{C(N)}X_T+CY_TX_T.\non
\end{align}
This gives by (\ref{equ:U0}) and (\ref{equ:U1estimateq}) that
\begin{equation}\label{equ:XT}
X_T\le C\exp\Big(\f {CT2^{2N}} {C(N)}+CY_T\Big)\Big(\f
{1+T2^{N(\f3q-\f3p+2)}} {C(N)}(1+X_T)+Y_T+X_T Y_T\Big).
\end{equation}

Due to
$\max\big(\frac{2}{\widetilde{p}}-\frac1q-\frac1p,\,\,\frac3{5q}-\frac3{5p}\big)<\epsilon$,
we can take $N$ big enough such that
\begin{align*}
\frac{CT^\f122^{N(\frac3{\widetilde{p}}-\frac3p+1)}}{C(N)}\ll 1,\quad
\frac{CT2^{N(\frac3q-\frac3p+2)}}{C(N)}\ll 1
\end{align*}
for any $T\le T_0$. By making use of the continuation argument, we
deduce from (\ref{equ:YT}) and (\ref{equ:XT}) that for any $T\le
T_0$,
\begin{align}\label{equ:XY-est}
X_T\le \f {C\big(1+T2^{N(\f3q-\f3p+2)}\big)} {C(N)},\quad
Y_T\le \frac {CT2^{N(\frac3q-\frac3p+2)}}{C(N)}.
\end{align}

Summing up (\ref{equ:U0}), (\ref{eq:U1-final}) and (\ref{equ:XY-est}), we conclude that
\begin{align*}
\|u(T_0)\|_{\dot B^{\frac3p-1}_{p,1}}&\ge\|U_1(T_0)\|_{\dot B^{\frac3p-1}_{p,1}}-\|U_0(T_0)\|_{\dot B^{\frac3p-1}_{p,1}}
-\|U_2(T_0)\|_{\dot B^{\frac3p-1}_{p,1}}\nonumber\\&\ge c2^{N(1-\frac3q-\frac3p-3\epsilon)}-CN2^{-\frac N2(\frac3q-\frac3p+\epsilon)}-C2^{N(\frac3{2q}-\frac3{2p}-\frac52\epsilon)}
\nonumber\\&\ge\frac c22^{N(1-\frac3q-\frac3p-3\epsilon)},
\end{align*}
if $N$ is taken sufficiently large, as
$1-\frac3q-\frac3p-3\epsilon>0$.  This completes the proof of
Theorem \ref{thm:illposed-bNS}.\ef

\section{Ill-posedness of the heat-conductive flows}

\subsection{Reformulation of the equation} We denote
$$a=\frac{\rho}{\bar\rho}-1,\quad u=U_0+\widetilde{U},$$
where $U_0$ is defined as  in \eqref{equ:U0}.
Then the system \eqref{equ:cNS} can be rewritten as
\begin{equation}\label{equ:heatflow-new}
\left\{
\begin{aligned}{}
&\p_ta+u\cdot\na a+\dv u(1+a)=0,\\
&\p_t \widetilde{U}-\mathcal{A}\widetilde{U}=-u\cdot\na u-L(a)\mathcal{A}u-R\nabla \theta-R\na a\theta+R\na a L(a)\theta,\\
&\p_t\theta+u\cdot\na\theta-\tilde{\kappa}\Delta\theta+\tilde{\kappa}L(a)\Delta\theta+\tilde{R}\theta\dv u=\frac{\tilde{\mu}}{2(a+1)}|\na u+(\na u)^\top|^2+\frac{\tilde{\lambda}}{a+1}|\dv u|^2,\\
&(a,\,\widetilde{U},\,\theta)|_{t=0}=\big(a_0,\,0,\,\theta_0\big),
\end{aligned}
\right.
\end{equation}
where $\tilde{\kappa}=\frac{\kappa}{c_V\bar{\rho}}$,
$\tilde{R}=\frac{R}{c_V}$, $\tilde{\mu}=\frac{\mu}{c_V\bar{\rho}}$,
$\tilde{\lambda}=\frac{\lambda}{c_V\bar{\rho}}$, and \beno
\mathcal{A}=\bar{\mu}\Delta+(\bar{\lambda}+\bar{\mu})\na \dv,\quad
L(a)=\frac{a}{1+a}. \eeno We   decompose $\theta$ into
\begin{align}\label{equ:decompositiontheta}
\theta(x,t)=\Theta_0+\theta_1+\theta_2,
\end{align}
where \begin{align*}
\Theta_0=&e^{\tilde{\kappa}\Delta t}\theta,
\\\theta_1=&\int_0^te^{\tilde{\kappa}\Delta (t-\tau)}\Big(\frac{\tilde{\mu}}{2}|\na U_0+(\na U_0)^\top|^2+
{\tilde{\lambda}}|\dv U_0|^2\Big)d\tau,\\
\theta_2=&\int_0^te^{\tilde{\kappa}\Delta
(t-\tau)}\Big(\frac{\tilde{\mu}}{2}|\na\widetilde{U}+(\na
\widetilde{U})^\top|^2+\tilde{\mu}(\na U_0+ (\na
U_0)^\top):(\na\widetilde{U}+(\na
\widetilde{U})^\top)\nonumber\\&\quad+
\tilde{\lambda}\big(|\dv\widetilde{U}|^2+2\dv
U_0\dv\widetilde{U}\big)-u\cdot\na\theta-\tilde{\kappa}L(a)\Delta\theta\nonumber\\&
\quad-\tilde{R}\theta\dv u-\frac{\tilde{\mu}}{2}L(a)|\na u+(\na
u)^\top|^2 -\tilde{\lambda}L(a)|\dv u|^2 \Big)d\tau.
\end{align*}

\subsection{The choice of initial data}

Let $\phi$ be as in \eqref{equ:phi} and $N\in \N$ be determined later. The initial velocity is chosen as
\begin{align*}
\widehat{u}_0(\xi)=\frac{1}{C(N)}\sum_{k=10}^N 2^{k(1-\frac3p)}\big(\phi(\xi-2^k{e_1})+\phi(\xi+2^k{e_1}),\,\, 0,
\,\,0\big),
\end{align*}
and the initial density and temperature is chosen as
\begin{align}
a_0=\frac{\rho_0}{\bar\rho}-1=\frac{1} {2^NC(N)}\mathcal{F}^{-1}(\phi)(x),\quad
\theta_0=\frac{1}{C(N)}\mathcal{F}^{-1}(\phi)(x).\non
\end{align}
Here $e_1=(1,0,0)$ and $C(N)=2^{\frac N2(\frac3q-\frac3p+\epsilon)}$.

\begin{Lemma}Let $p>3$. There exist $\epsilon>0$ and $(\widetilde{p}, q)$ satisfying
\begin{align*}
&2<q<3,\quad 3<\widetilde{p}<p,\quad \frac3{\widetilde{p}}+\frac3q-2>0,\\
&\max\Big\{\frac{2}{\widetilde{p}}-\frac1q-\frac1p,\,\frac35\Big(\frac1{q}-\frac1{p}\Big)\Big\}<\epsilon<\frac23-\frac1q-\frac1p.
\end{align*}
\end{Lemma}

Throughout this section, we will fix such a triplet $(\epsilon,
\widetilde{p}, q)$. It is easy to check that
\begin{align}\label{equ:u0-h}
\|a_0\|_{\dot B^{\frac3q}_{q,1}}+\|\theta_0\|_{\dot
B^{-2+\frac3q}_{q,1}}\le \f C{C(N)},\quad \|u_0\|_{\dot
B^{-1+\frac3p}_{p,1}}\le \frac {CN} {C(N)}.
\end{align}
Furthermore, it holds that
\begin{align}\label{equ:U0-h}
\|U_0\|_{\widetilde{L}^\rho_T\dot B^{\sigma}_{r,1}}\le
CT^{\frac1{\rho_1}}\frac{2^{N(\sigma-\frac3p+1-\frac2{\rho_2})}}{C(N)}
\end{align}
for any $r,\rho, \rho_1, \rho_2\in [1,\infty]$ and
$\sigma>\frac3p-1+\frac2{\rho_2}$ with
$\frac1\rho=\frac1{\rho_1}+\frac1{\rho_2}$.

\subsection{The lower bound estimate of $\|\theta_1\|_{\dot B^{\frac3p-2}_{p,1}}$}

By the same argument as the one used  to derive \eqref{equ:U1_p} and
\eqref{equ:U1>}, we have
\begin{align}\label{equ:theta1}
&\|\theta_1(t)\|_{\dot B^{\frac3p-2}_{p,1}}\ge c\Big|\int_{\R^3}\varphi(2^4\xi) \widehat{\theta}_1(t,\xi)d\xi\Big|
\ge c|\theta_{11}+\theta_{12}|,
\end{align}
where $\theta_{11}$ and $\theta_{12}$ are given by
\begin{align*}
&\mathfrak{\theta}_{11}=\tilde{\lambda}\int_{}\!\int_0^t{\varphi(2^4\xi)} e^{-{\tilde\kappa(t-\tau)|\xi|^2}}
\mathcal{F}\big(|\dv U_0|^2\big)(\tau,\xi)d\tau d\xi,
\\&\mathbf{\theta}_{12}=\frac{\tilde{\mu}}{2}\int_{}\!\int_0^t\varphi(2^4\xi) e^{-{\tilde\kappa(t-\tau)|\xi|^2}}\mathcal{F}\big(|\na U_0+(\na U_0)^\top|^2\big)(\tau,\xi)d\tau d\xi.
\end{align*}

$\bullet$\,\,{\bf The estimate of $\theta_{11}$} \vspace{0.1cm}

Recalling the definition of $U_0$, we get
\begin{align*}
\theta_{11}=&-\widetilde{\lambda}\int_{}\!\!\int_{}\!\!\int_0^t\varphi(2^4\xi)
e^{-{\tilde\kappa t|\xi|^2}}
e^{\tilde\kappa\tau|\xi|^2-{\bar\nu\tau|\xi-\eta|^2}-{\bar\nu\tau|\eta|^2}}
(\xi-\eta)_j\widehat{u_0}^j(\xi-\eta)\eta_\ell\widehat{u_0}^\ell(\eta)d\eta d\tau d\xi\\
=&\frac{\tilde{\lambda}}{C(N)^2}\sum_{k=10}^N
2^{2(1-\frac3p)k}\int_{}\!\!\int\varphi(2^4\xi)\frac{e^{-\bar{\kappa}|\xi|^2t}-e^{(-\bar{\nu}|\xi-\eta|^2-\bar{\nu}|\eta|^2)t}}
{\bar{\kappa}|\xi|^2-\bar{\nu}|\xi-\eta|^2-\bar{\nu}|\eta|^2}
({\xi}_1-\eta_1)\eta_1\nonumber\\
&\quad\times\big(\phi(\xi-\eta+2^ke_1)\phi
(\eta-2^ke_1)+\phi(\xi-\eta-2^ke_1)\phi (\eta+2^ke_1)\big)d\eta d\xi.
\end{align*}
Making a change of variable, we find that
\begin{align*}
\theta_{11}=\frac{\tilde{\lambda}}{C(N)^2}\sum_{k=10}^N
2^{2(1-\frac3p)k}\int_{}\!\!\int&\varphi(2^4\xi)\frac{e^{-\bar{\kappa}|\xi|^2t}-e^{(-\bar{\nu}
|\xi-\eta-2^ke_1|^2-\bar{\nu}|\eta+2^ke_1|^2)t}}{\bar{\kappa}|\xi|^2
-\bar{\nu}|\xi-\eta-2^ke_1|^2-\bar{\nu}|\eta+2^ke_1|^2}\nonumber\\
&\times
({\xi}_1-\eta_1-2^k)(\eta_1+2^k)\phi(\xi-\eta)\phi(\eta)d\eta d\xi\\
+\frac{\tilde{\lambda}}{C(N)^2}\sum_{k=10}^N
2^{2(1-\frac3p)k}\int_{}\!\!\int&\varphi(2^4\xi)\frac{e^{-\bar{\kappa}|\xi|^2t}-e^{(-\bar{\nu}
|\xi-\eta+2^ke_1|^2-\bar{\nu}|\eta-2^ke_1|^2)t}}{\bar{\kappa}|\xi|^2
-\bar{\nu}|\xi-\eta+2^ke_1|^2-\bar{\nu}|\eta-2^ke_1|^2}\nonumber\\
&\times({\xi}_1-\eta_1+2^k)(\eta_1-2^k)\phi(\xi-\eta)\phi(\eta)d\eta d\xi.
\end{align*}
Using  Taylor's formula, we deduce
\begin{align}\label{equ:theta11}
\theta_{11}(t)=\frac{2\tilde{\lambda}}{C(N)^2}\sum_{k=10}^N
2^{2(1-\frac3p)k}\int_{}\!\!\int&\varphi(2^4\xi)\big(-t+{O(2^{2k}t^2)}\big)\nonumber\\
&\times(-2^{2k}+O(2^k))\phi(\xi-\eta)\phi(\eta)d\eta d\xi.
\end{align}

$\bullet$\,\,{\bf The estimate of $\theta_{12}$} \vspace{0.1cm}

Notice that
\beno
\mathcal{F}\big(\na U_0+(\na U_0)^\top\big)(\eta)
=-\big(\widehat{H}_{\ell m}^1(\eta)+\widehat{H}_{\ell m}^2(\eta)+\widehat{H}_{\ell m}^3(\eta)\big),
\eeno
where $\widehat{H}_{\ell m}^i(i=1,2,3)$ is given by
\begin{align*}
\widehat{H}_{\ell m}^1(\eta)=&\frac{2i}{|\eta|^2}e^{-\bar\nu|\eta|^2\tau}\eta_\ell\eta_m\eta_j\widehat{u}_0^j(\eta),\\
\widehat{H}_{\ell m}^2(\eta)=&-\frac{2i}{|\eta|^2}e^{-\bar\mu|\eta|^2\tau}\eta_\ell\eta_m\eta_j\widehat{u}_0^j(\eta),\\
\widehat{H}_{\ell
m}^3(\eta)=&\frac{i}{|\eta|^2}e^{-\bar\mu|\eta|^2\tau}\eta_j^2(\eta_\ell\widehat{u}_0^m(\eta)
+\eta_m\widehat{u}_0^\ell(\eta)),
\end{align*}
we   write $\theta_{12}$ as
\begin{align}
\theta_{12}&=\frac{\tilde{\mu}}{2}\sum_{J, J'=1}^3\int_{}\!\int_{}\!\int_0^t\varphi(2^4\xi) e^{-{\tilde\kappa(t-\tau)|\xi|^2}}
\widehat{H}_{\ell m}^J(\xi-\eta)\widehat{H}_{\ell m}^{J'}(\eta)d\tau d\eta d\xi\nonumber\\
&\triangleq \sum_{J, J'=1}^3\mathcal{H}^{JJ'}.\non
\end{align}

First of all, we consider $\mathcal{H}^{11}$.
\begin{align*}\label{}
\mathcal{H}^{11}(t)=-2\tilde{\mu}\int_{}\!\!\int_{}\!\!\int_0^t&\varphi(2^4\xi)\frac{e^{-{\tilde\kappa(t-\tau)|\xi|^2}}}
{|\xi-\eta|^2|\eta|^2}e^{-{\bar\nu\tau|\xi-\eta|^2}-{\bar\nu\tau|\eta|^2}}
(\xi-\eta)_m(\xi-\eta)_\ell(\xi-\eta)_j
\nonumber\\&\times\eta_m\eta_\ell\eta_{j'}\widehat{u}_0^j(\xi-\eta)\widehat{u}_0^{j'}(\eta)d\tau d\eta d\xi.
\end{align*}
Making a change of variable and computing the above integral with
respect to the time, we deduce
\begin{align}
\mathcal{H}^{11}(t)=&\frac{2\tilde{\mu}}{C(N)^2}\sum_{k=10}^N
2^{2(1-\frac3p)k}\int_{}\!\!\int\varphi(2^4\xi)
\frac{e^{-\tilde{\kappa}|\xi|^2t}-e^{(-\bar{\nu}
|\xi-\eta-2^ke_1|^2-\bar{\nu}|\eta+2^ke_1|^2)t}}{\tilde{\kappa}|\xi|^2
-\bar{\nu}|\xi-\eta-2^ke_1|^2-\bar{\nu}|\eta+2^ke_1|^2}
\nonumber\\&\quad\times\frac{(\xi-\eta-2^ke_1)_{\ell}}{|\xi-\eta-2^ke_1|^2}\frac{(\eta+2^ke_1)_{\ell}}{|\eta+2^ke_1|^2}
(\xi-\eta-2^ke_1)_{m}(\eta+2^ke_1)_{m} \nonumber\\&\quad\times
(\xi-\eta-2^ke_1)_{1}(\eta+2^ke_1)_{1}\phi(\xi-\eta)\phi(\eta)d\eta d\xi+\textrm{Similar term}.\non
\end{align}
Using Taylor's formula, we infer
\begin{align}
\mathcal{H}^{11}(t)=\frac{4\tilde{\mu}}{C(N)^2}\sum_{k=10}^N
2^{2(1-\frac3p)k}\int_{}\!\!\int&\varphi(2^4\xi)\big(-t+{O(2^{2k}t^2)}\big)\nonumber\\
&\times(-2^{2k}+O(2^{k}))\phi(\xi-\eta)\phi(\eta)d\eta d\xi.\label{equ:H11}
\end{align}

By a direct calculation, we find
\begin{align*}
&\mathcal{H}^{33}(t)+\mathcal{H}^{32}(t)=-2\tilde{\mu}\int_{}\!\!\int_{}\!\!\int_0^t\varphi(2^4\xi)
\frac{e^{-{\tilde\kappa(t-\tau)|\xi|^2}}}
{|\xi-\eta|^2|\eta|^2}e^{-{\bar\mu\tau|\xi-\eta|^2}-{\bar\mu\tau|\eta|^2}}(\xi-\eta)_j^2\nonumber\\&\times
\Big\{\!(\eta_2^2+\eta_3^2)\Big(2(\xi-\eta)_1\eta_1
+(\xi-\eta)_2\eta_2+(\xi-\eta)_3\eta_3\Big)
-\eta_1^2\Big((\xi-\eta)_2\eta_2+(\xi-\eta)_3\eta_3\Big)
\Big\}\\
&\qquad\qquad\times\hat{u}_0^1(\xi-\eta)\hat{u}_0^{1}(\eta)d\tau d\eta d\xi..
\end{align*}
Integrating on the time and making a change of variable, we deduce
\begin{align*}
\mathcal{H}&^{33}(t)+\mathcal{H}^{32}(t)\nonumber\\=&\frac2{C(N)^2}\sum_{k=10}^N
2^{2(1-\frac3p)k}\int_{}\!\!\int\varphi(2^4\xi)
\frac{e^{-\tilde{\kappa}|\xi|^2t}-e^{(-\bar{\mu}
|\xi-\eta-2^ke_1|^2-\bar{\mu}|\eta+2^ke_1|^2)t}}{\tilde{\kappa}|\xi|^2
-\bar{\mu}|\xi-\eta-2^ke_1|^2-\bar{\mu}|\eta+2^ke_1|^2}
\nonumber\\&\quad\times\frac{(\xi-\eta-2^ke_1)_{j}^2}{|\xi-\eta-2^ke_1|^2|\eta+2^ke_1|^2}
O(2^{2k})\phi(\xi-\eta)\phi(\eta)d\eta d\xi+\textrm{Simliar term},
\end{align*}
from which, it follows that
\begin{align}
|\mathcal{H}^{33}(t)+\mathcal{H}^{32}(t)|\le \frac C{C(N)^2}\sum_{k=10}^N2^{-\frac6pk}\le \frac C {C(N)^2}.
\end{align}
Furthermore, we also have
\begin{align*}
\mathcal{H}&^{22}(t)+\mathcal{H}^{23}(t)\nonumber\\=&-2\widetilde{\mu}\int_{}\!\!\int_{}\!\!\int_0^t\varphi(2^4\xi)
\frac{e^{-{\tilde\kappa(t-\tau)|\xi|^2}}}
{|\xi-\eta|^2|\eta|^2}e^{-{\bar\mu\tau|\xi-\eta|^2}-{\bar\mu\tau|\eta|^2}}(\xi-\eta)_\ell(\xi-\eta)_1\eta_\ell
\nonumber\\
&\qquad\times\big\{(\xi-\eta)_2\eta_2\eta_1+(\xi-\eta)_3\eta_3\eta_1-(\xi-\eta)_1(\eta_2^2+\eta_3^2)\big\}\widehat{u}_0^1(\xi-\eta)\widehat{u}_0^1(\eta)d\tau d\eta d\xi,
\end{align*}
and a similar representation for $\mathcal{H}^{12}(t)+\mathcal{H}^{13}(t)$, hence,
\begin{align}
|\mathcal{H}^{22}(t)+\mathcal{H}^{23}(t)|+|\mathcal{H}^{12}(t)+\mathcal{H}^{13}(t)|\le
\frac C{C(N)^2}\sum_{k=10}^N2^{-\frac6pk}\le \frac C {C(N)^2}.
\end{align}
For $\mathcal{H}^{21}(t)+\mathcal{H}^{31}(t)$, we have
\begin{align*}
\mathcal{H}&^{21}(t)+\mathcal{H}^{31}(t)\nonumber\\=&2\widetilde{\mu}\int_{}\!\!\int_{}\!\!\int_0^t\varphi(2^4\xi)
\frac{e^{-{\tilde\kappa(t-\tau)|\xi|^2}}}
{|\xi-\eta|^2|\eta|^2}e^{-{\bar\mu\tau|\xi-\eta|^2}-{\bar\nu\tau|\eta|^2}}(\xi-\eta)_\ell\eta_\ell\eta_1
\nonumber\\&\times\big\{((\xi-\eta)_2\eta_2+(\xi-\eta)_3\eta_3)(\xi-\eta)_1
-((\xi-\eta)^2_2+(\xi-\eta)^2_3)\eta_1\big\}\widehat{u}_0^1(\xi-\eta)\widehat{u}_0^1(\eta)d\tau d\eta d\xi,
\end{align*}
hence,
\begin{align}
|\mathcal{H}^{21}(t)+\mathcal{H}^{31}(t)|\le\frac C {C(N)^2}.\label{equ:H21+31}
\end{align}

Summing up (\ref{equ:theta1})--(\ref{equ:H21+31}), we conclude that
\begin{align} \|\theta_1(t)\|_{\dot B^{\frac3p-2}_{p,1}}\ge&
|\theta_{11}+\mathcal{H}^{11}(t)|-
{\Big|\sum_{1\le J, J'\le3,\,JJ'\neq11}\mathcal{H}^{JJ'}(t)\Big|}\non\\
\ge& \f{c2^{2(2-\f3p)N}t}{C(N)^2}-\f
{C\big(1+2^{6(1-\f1p)N}t^2+2^{(3-\f6p)N}t\big)}{C(N)^2}.\non
\end{align} Especially, due to $2-\frac3q-\frac3p>3\epsilon$,  take
$t\sim 2^{-2(1+\epsilon)N}$  to obtain
\begin{align}\label{equ:theta1final}
\|\theta_1(t)\|_{\dot B^{\frac3p-2}_{p,1}}\ge c2^{N(2-\frac3q-\frac3p-3\epsilon)}
\end{align}
for some $c>0$ independent of $N$.

\subsection{The estimate of $\|\theta_1\|_{L^1_T\dot B^{\frac3q}_{q,1}\cap \widetilde{L}^2_T\dot B^{\frac3q-1}_{q,1}}$}

Since $q<3$, we get by Lemma \ref{Lem:binestis>0} that
\begin{align*}
\|(\na U_0)^2\|_{L^1_T\dot B^{\frac3q-1}_{q,1}}&\le C(\|\nabla
U_0\|_{L^\infty_TL^\infty}\|\nabla U_0\|_{L^1_T\dot
B^{\frac3q-1}_{q,1}}+\|\nabla U_0\|_{L^1_T\dot
B^{\frac3q-1}_{q,1}}\|\nabla U_0\|_{L^\infty_TL^\infty})
\nonumber\\&\le CT^\frac12\|U_0\|_{L^\infty_T\dot
B^{\frac3p+1}_{p,1}}\|U_0\|_{\widetilde{L}^2_T\dot
B^{\frac3q}_{q,1}}.
\end{align*}
By H\"{o}lder's inequality and Proposition \ref{Prop:heatflow}, we
get
\begin{align}\label{equ:theta1estimate0}
\|\theta_1\|_{L^1_T\dot
B^{\frac3q}_{q,1}}+\|\theta_1\|_{\widetilde{L}^2_T\dot
B^{\frac3q-1}_{q,1}}&\le C\big(
T^{\frac12}\|\theta_1\|_{\widetilde{L}^2_T\dot
B^{\frac3q}_{q,1}}+\|\theta_1\|_{\widetilde{L}^2_T\dot
B^{\frac3q-1}_{q,1}}\big)
\nonumber\\&\le CT^{\frac12}\|(\na U_0)^2\|_{L^1_T\dot B^{\frac3q-1}_{q,1}}\non\\
&\le  CT\|U_0\|_{L^\infty_T\dot
B^{\frac3p+1}_{p,1}}\|U_0\|_{\widetilde{L}^2_T\dot
B^{\frac3q}_{q,1}}.\non
\end{align}
This along with \eqref{equ:U0-h} gives \ben\label{equ:theta1-f2}
\|\theta_1\|_{L^1_T\dot B^{\frac3q}_{q,1}}+\|\theta_1\|_{L^2_T\dot
B^{\frac3q-1}_{q,1}} \le CT\frac{2^{N(\frac3q-\frac3p+2)}}{C(N)^2}.
\een

\subsection{Proof of Theorem \ref{thm:illposed-heat}}

We denote
\beno
&&X_T=\|{a}\|_{\widetilde{L}^\infty_T\dot B^{\frac3q}_{q,1}},\quad Y_T=\|\widetilde{U}\|_{\widetilde{L}^\infty_T\dot B^{\frac3q-1}_{q,1}}+
\|\widetilde{U}\|_{L^1_T\dot B^{\frac3q+1}_{q,1}},\\
&&Z_T=\|\theta_2\|_{\widetilde{L}^\infty_T\dot B^{\frac3q-2}_{q,1}}+\|\theta_2\|_{L^1_T\dot B^{\frac3q}_{q,1}}.
\eeno

{\bf Step 1.} The estimate of $X_T$\vspace{0.1cm}

It follows from Proposition \ref{Prop:transport}, Lemma \ref{Lem:binestis>0} and (\ref{equ:U0-h}) that
\ben
&&X_T\le C \exp(\|\na u\|_{L^1_T\dot B^{\frac3p}_{p,1}})
\big(\|{a}_0\|_{\dot B^{\frac3q}_{q,1}}+\|\dv(U_0+\widetilde{U})+
\dv (U_0+\widetilde{U})\,{a}\|_{L^1_T\dot B^{\frac3q}_{q,1}}\big)\nonumber\\
&&\le
C\exp\Big(\frac{T2^{2N}}{C(N)}+Y_T\Big)\Big(\frac{1}{C(N)}+\frac{T2^{N(\frac{3}{q}-\frac3p+2)}}{C(N)}(1+X_T)+Y_T
+Y_TX_T\Big).\label{equ:XT-h} \een Similarly, we have
\begin{align}
\|a\|_{{\widetilde{L}}^\infty_T\dot B^{\frac3q-1}_{q,1}}\le& C\exp\Big(\frac{T2^{2N}}{C(N)}+Y_T\Big)\nonumber\\
&\times\Big(\frac{1}{2^NC(N)}+\frac{T2^{N(\frac{3}{q}-\frac3p+1)}}{C(N)}(1+X_T)+TY_T+T^\f12Y_TX_T\Big).
\end{align}

{\bf Step 2.} The estimate of $Y_T$\vspace{0.1cm}

We infer from Proposition \ref{Prop:heatflow} that
\begin{align}\label{equ:YT-h1}
Y_T\le C\|u\cdot\na u+L(a)\mathcal{A}u+R\nabla \theta+R\na a\theta-R\na a L(a)\theta\|_{L^1_T\dot B^{\frac3q-1}_{q,1}}.
\end{align}
By Lemma \ref{Lem:binesti} and (\ref{equ:U0-h}), we get
\begin{align}\label{equ:YT-h2}
\|L(a)\mathcal{A} u\|_{L^1_T\dot B^{\frac3q-1}_{q,1}}&\le C\|a\|_{L^\infty_T\dot B^{\frac3q}_{q,1}}
\|\na^2 (U_0+\widetilde{U})\|_{L^1_T\dot B^{\frac3q-1}_{q,1}}\nonumber\\
&\le CX_T\Big(\frac{T2^{N(\frac{3}{q}-\frac3p+2)}}{C(N)}+Y_T\Big).
\end{align}
Using Lemma \ref{Lem:nonesti}, (\ref{equ:u0-h}) and
(\ref{equ:theta1-f2}) yields that
\begin{align}\label{equ:YT-h3}
&\|\na a\theta\|_{L^1_T\dot B^{\frac3q-1}_{q,1}}+\|\na a L(a)\theta\|_{L^1_T\dot B^{\frac3q-1}_{q,1}}\nonumber\\
&\le C\|\na a\|_{L^\infty_T\dot B^{\frac3q-1}_{q,1}}\|\theta\|_{L^1_T\dot B^{\frac3q}_{q,1}}+C\|\na a\|_{L^\infty_T\dot B^{\frac3q-1}_{q,1}}\|L(a)\|_{L^\infty_T\dot B^{\frac3q}_{q,1}}\|\theta\|_{L^1_T\dot B^{\frac3q}_{q,1}}\nonumber\\&
\le CX_T(1+X_T)^3\Big(\f 1{C(N)}+\frac{T2^{N(\frac{3}{q}-\frac3p+2)}}{C(N)^2}+Z_T\Big).
\end{align}
On the other hand, by (\ref{Bonydecom}), Lemma \ref{lem:biest-n} and (\ref{equ:U0-h}), we get
\begin{align}
\|U_0\cdot \nabla U_0\|_{L^1_T\dot B^{\frac3q-1}_{q,1}}
&\le C\|U_0\|_{\widetilde{L}^2_T\dot B^{\frac3p}_{p,1}}\|U_0\|_{\widetilde{L}^2_T\dot B^{\frac3q}_{q,1}}
\le C\frac{T2^{N(\frac{3}{q}-\frac3p+2)}}{C(N)^2},
\end{align}
and
\begin{align}
\|U_0\cdot \nabla \widetilde{U}\|_{L^1_T\dot B^{\frac3q-1}_{q,1}}&\le C\|U_0\|_{\widetilde{L}^2_T\dot B^{\frac3p}_{p,1}}\|\nabla \widetilde{U}\|_{\widetilde{L}^2_T\dot B^{\frac3q-1}_{q,1}}
+C\|U_0\|_{L^1_T\dot B^{\frac3p+1}_{p,1}}\|\nabla \widetilde{U}\|_{\widetilde{L}^\infty_T\dot B^{\frac3q-2}_{q,1}}\nonumber\\
&\le C\frac{T^{\frac12}2^N+T2^{2N}}{C(N)}Y_T.
\end{align}
Thanks to Lemma \ref{Lem:binesti}, we have
\begin{align}\label{equ:YT-h6}
\|\widetilde{U}\cdot \na u\|_{L^1_T\dot B^{\frac3q-1}_{q,1}}\le& C\|\widetilde{U}\|_{L^\infty_T\dot B^{\frac3q-1}_{q,1}}\big(
\|\na U_0\|_{L^1_T\dot B^{\frac3q}_{q,1}}+\|\na \widetilde{U}\|_{L^1_T\dot B^{\frac3q}_{q,1}}\big)\nonumber\\
\le& CY_T\Big(\frac{T2^{N(\frac{3}{q}-\frac3p+2)}}{C(N)}+Y_T\Big).
\end{align}

Summing up (\ref{equ:YT-h1})--(\ref{equ:YT-h6}), we conclude that
\begin{align}\label{equ:YT-h-f}
Y_T\le& \frac{CT(1+2^{N(\frac{3}{q}-\frac3p+2)}/C(N))}{C(N)}+\f {C(1+T2^{N(\f3q-\f3p+2)})} {C(N)}(1+X_T)^3X_T+Z_T\non\\
&+\f {C\big(T^\f122^N+T2^{N(\f3q-\f3p+2)}\big)}
{C(N)}Y_T+C(1+X_T)^3X_TZ_T+CX_TY_T+CY_T^2.
\end{align}

{\bf Step 3.} The estimate of $Z_T$\vspace{0.1cm}

It follows from Proposition \ref{Prop:heatflow} that
\begin{align}\label{equ:ZT-1}
Z_T\le& C\Big\{\|(\na \widetilde{U})^2\|_{L^1_T\dot B^{\frac3q-2}_{q,1}}+\|\na U_0\na \widetilde{U}\|_{L^1_T\dot B^{\frac3q-2}_{q,1}}+\|u\cdot\na \theta\|_{L^1_T\dot B^{\frac3q-2}_{q,1}}\nonumber\\&\quad+\|L(a)\Delta\theta\|_{L^1_T\dot B^{\frac3q-2}_{q,1}}+\|\theta\dv u\|_{L^1_T\dot B^{\frac3q-2}_{q,1}}+\|L(a)(\na U_0)^2\|_{L^1_T\dot B^{\frac3q-2}_{q,1}}\nonumber\\&\quad+\|L(a)(\na \widetilde{U})^2\|_{L^1_T\dot B^{\frac3q-2}_{q,1}}
+\|L(a)\na U_0\na \widetilde{U}|\|_{L^1_T\dot B^{\frac3q-2}_{q,1}}\Big\}.
\end{align}

Next we estimate each term on the right hand side of (\ref{equ:ZT-1}).
Thanks to Lemma \ref{Lem:binesti} and (\ref{equ:U0-h}), we get
\begin{align}\label{equ:ZT-2}
&\|(\na \widetilde{U})^2\|_{L^1_T\dot B^{\frac3q-2}_{q,1}}+\|\na U_0\na \widetilde{U}\|_{L^1_T\dot B^{\frac3q-2}_{q,1}}\nonumber\\
&\le C\|\na \widetilde{U}\|_{\widetilde{L}^2_T\dot B^{\frac3q-1}_{q,1}}^2+C\|\na U_0\|_{L^1_T\dot B^{\frac3q}_{q,1}}\|\na \widetilde{U}\|_{L^\infty_T\dot B^{\frac3q-2}_{q,1}}\nonumber\\
&\le CY_T^2
+\frac{CT2^{N(\frac{3}{q}-\frac3p+2)}}{C(N)}Y_T.
\end{align}
Due to $\f 3p-\f3q+2>0$ and $\f 3 {\widetilde{p}}+\f3q-2>0$,  we
apply (\ref{Bonydecom}), Lemma \ref{lem:biest-n}, (\ref{equ:U0-h})
and (\ref{equ:theta1-f2}) to obtain
\begin{align}\label{equ:ZT-3}
&\|U_0\cdot\na \theta\|_{L^1_T\dot B^{\frac3q-2}_{q,1}}\non\\
&\le\|U_0\|_{\widetilde{L}^2_T\dot B^{\frac3p}_{p,1}}\|\na
\theta\|_{\widetilde{L}^2_T\dot
B^{\frac3q-2}_{q,1}}+\|U_0\|_{\widetilde{L}^2_T\dot
B^{\frac3{\widetilde{p}}}_{\widetilde{p},1}}\|\na
\theta\|_{\widetilde{L}^2_T\dot B^{\frac3q-2}_{q,1}}\nonumber\\&\le
CT^{\frac12}\Big(\frac{2^N}{C(N)}+\frac{2^{N(\frac{3}{\widetilde{p}}-\frac3p+1)}}{C(N)}\Big)
\Big(\|\Theta_0\|_{\widetilde{L}^2_T\dot
B^{\frac3q-1}_{q,1}}+\|\theta_1\|_{\widetilde{L}^2_T\dot
B^{\frac3q-1}_{q,1}}+\|\theta_2\|_{\widetilde{L}^2_T\dot
B^{\frac3q-1}_{q,1}}\Big) \nonumber\\&\le
CT^{\frac12}\frac{2^{N(\frac{3}{\widetilde{p}}-\frac3p+1)}}{C(N)}
\Big(\f 1
{C(N)}+\frac{T2^{N(\frac{3}{q}-\frac3p+2)}}{C(N)}+Z_T\Big).
\end{align}
On the other hand, we have by Lemma \ref{Lem:binesti}
\begin{align}
\|\widetilde{U}\cdot\na \theta\|_{L^1_T\dot B^{\frac3q-2}_{q,1}}&\le C\|\widetilde{U}\|_{\widetilde{L}^2_T\dot B^{\frac3q}_{q,1}}\|\na \theta\|_{\widetilde{L}^2_T\dot B^{\frac3q-2}_{q,1}}\nonumber\\
&\le CY_T\Big(\f 1
{C(N)}+\frac{T2^{N(\frac{3}{q}-\frac3p+2)}}{C(N)}+Z_T\Big),\label{equ:ZT-4}
\end{align}
and
\begin{align}
\|L(a)\Delta\theta\|_{L^1_T\dot B^{\frac3q-2}_{q,1}}&\le C\|L(a)\|_{L^\infty_T\dot B^{\frac3q}_{q,1}}\|\Delta \theta\|_{L^1_T\dot B^{\frac3q-2}_{q,1}}\non\\
&\le C(1+X_T)^3X_T\Big(\f 1 {C(N)}+\frac{T2^{N(\frac{3}{q}-\frac3p+2)}}{C(N)}+Z_T\Big),\label{equ:ZT-5}
\end{align}
and
\begin{align}
\|\theta\dv \widetilde{U}\|_{L^1_T\dot B^{\frac3q-2}_{q,1}}&\le C\|\theta\|_{L^1_T\dot B^{\frac3q}_{q,1}}\|\dv \widetilde{U}\|_{L^\infty_T\dot B^{\frac3q-2}_{q,1}}\nonumber
\\&\le C\Big(\f 1 {C(N)}+\frac{T2^{N(\frac{3}{q}-\frac3p+2)}}{C(N)}+Z_T\Big)Y_T.
\end{align}
Collecting (\ref{Bonydecom}), (\ref{equ:U0-h})  and
(\ref{equ:theta1-f2}) with Lemma \ref{lem:biest-n} implies that
\begin{align}
&\|(\Theta_0+\theta_2)\dv U_0\|_{L^1_T\dot B^{\frac3q-2}_{q,1}}\nonumber\\
&\le C\|\dv U_0\|_{L^1_T\dot B^{\frac3p}_{p,1}}\|\Theta_0+\theta_2\|_{L^\infty_T\dot B^{\frac3q-2}_{q,1}}+C\|\dv U_0\|_{L^1_T\dot B^{\frac3{\widetilde{p}}}_{\widetilde{p},1}}\|\Theta_0+\theta_2\|_{L^\infty_T\dot B^{\frac3q-2}_{q,1}}\nonumber\\&\le CT\Big(\frac{2^{2N}}{C(N)}+\frac{2^{N(\frac{3}{\widetilde{p}}-\frac3p+2)}}{C(N)}\Big)
\Big(\|\Theta_0\|_{L^\infty_T\dot B^{\frac3q-2}_{q,1}}+\|\theta_2\|_{L^\infty_T\dot B^{\frac3q-2}_{q,1}}\Big)
\nonumber\\&\le CT\frac{2^{N(\frac{3}{\widetilde{p}}-\frac3p+2)}}{C(N)}
\big(\f 1 {C(N)}+Z_T\big), \label{equ:ZT-6}
\end{align}
and
\begin{align}
&\|\theta_1\dv U_0\|_{L^1_T\dot B^{\frac3q-2}_{q,1}}\nonumber\\
&\le C\|\dv U_0\|_{\widetilde{L}^2_T\dot B^{\frac3p-1}_{p,1}}\|\theta_1\|_{\widetilde{L}^2_T\dot B^{\frac3q-1}_{q,1}}+C\|\dv U_0\|_{\widetilde{L}^2_T\dot B^{\frac3{\widetilde{p}}-1}_{\widetilde{p},1}}\|\theta_1\|_{\widetilde{L}^2_T\dot B^{\frac3q-1}_{q,1}}\nonumber\\&\le CT^{\frac12}\Big(\frac{2^{N}}{C(N)}+\frac{2^{N(\frac{3}{\widetilde{p}}-\frac3p+1)}}{C(N)}\Big)
\frac{T2^{N(\frac{3}{q}-\frac3p+2)}}{C(N)^2},\label{equ:ZT-7}
\end{align}
and
\begin{align}
&\|L(a)(\na U_0)^2\|_{L^1_T\dot B^{\frac3q-2}_{q,1}}\nonumber\\&\le C\Big(\|\na U_0\|_{L^\infty_T\dot B^{\frac3p-1}_{p,1}}\|L(a)\na U_0\|_{L^1_T\dot B^{\frac3q-1}_{q,1}}+\|\na U_0\|_{L^1_T\dot B^{\frac3{\widetilde{p}}-1}_{\widetilde{p},1}}\|L(a)\na U_0\|_{L^\infty_T\dot B^{\frac3q-1}_{q,1}}\Big)
\nonumber\\&\le C\|\na U_0\|_{L^1_T\dot B^{\frac3{\widetilde{p}}-1}_{{\widetilde{p}},1}}\Big(\|\na U_0\|_{L^\infty_T\dot B^{\frac3p}_{p,1}}\|L(a)\|_{L^\infty_T\dot B^{\frac3q-1}_{q,1}}+\|\na U_0\|_{L^\infty_T\dot B^{\frac3{\widetilde{p}}}_{\widetilde{p},1}}\|L(a)\|_{L^\infty_T\dot B^{\frac3q-1}_{q,1}}\Big)
\nonumber\\&\le C\frac{T2^{N(\frac6{\widetilde{p}}-\frac6p+3)}}{C(N)^2}(1+X_T)^3\|a\|_{L^\infty_T\dot B^{\frac3q-1}_{q,1}}.
\end{align}
By Lemma \ref{Lem:binesti} and (\ref{equ:U0-h}), we get
\begin{align}\label{equ:ZT-9}
&\|L(a)\na U_0\na \widetilde{U}\|_{L^1_T\dot B^{\frac3q-2}_{q,1}}+\|L(a)(\na \widetilde{U})^2\|_{L^1_T\dot B^{\frac3q-2}_{q,1}}\nonumber\\&\le
C\|L(a)\|_{L^\infty_T\dot B^{\frac3q}_{q,1}}\big(\|\na U_0\na \widetilde{U}\|_{L^1_T\dot B^{\frac3q-2}_{q,1}}+\|(\na \widetilde{U})^2\|_{L^1_T\dot B^{\frac3q-2}_{q,1}}\big)\nonumber\\&\le
C\|a\|_{L^\infty_T\dot B^{\frac3q}_{q,1}}\big(\|\na U_0\|_{L^1_T\dot B^{\frac3q}_{q,1}}\|\na \widetilde{U}\|_{L^\infty_T\dot B^{\frac3q-2}_{q,1}}+\|\na \widetilde{U}\|_{L^1_T\dot B^{\frac3q}_{q,1}}\|\na \widetilde{U}\|_{L^\infty_T\dot B^{\frac3q-2}_{q,1}}\big)\nonumber\\
&\le
CX_T\Big(\frac{T2^{N(\frac{3}{q}-\frac3p+2)}}{C(N)}Y_T+Y_T^2\Big).
\end{align}

Summing up (\ref{equ:ZT-1})--(\ref{equ:ZT-9}), we deduce that
\begin{align}
Z_T\le& \f {C\big(T^\f122^{N(\f 3 {\widetilde{p}}-\f3p+1)}+T2^{N(\f
3 {\widetilde{p}}-\f3p+2)}
+T^\f322^{N(\f 3 {\widetilde{p}}+\f3q-\f6p+3)}\big)} {C(N)^2}\non\\
&+C\Big(Y_T+ \f {C\big(T^\f122^{N(\f 3 {\widetilde{p}}-\f3p+1)}+T2^{N(\f 3 {\widetilde{p}}-\f3p+2)}\big)} {C(N)} \Big)Z_T
+CY_T^2+CX_TY_T^2\non\\
&+C(1+X_T)^3\Big(Z_T+\f {1+T2^{N(\f3q-\f3p+2)}} {C(N)}+Y_T\f {T2^{N(\f3q-\f3p+2)}} {C(N)}\Big)(X_T+Y_T)\non\\
&+\f {CT2^{N(\f6{\widetilde{p}}-\f6p+3)}} {C(N)^2}(1+X_T)^3\|a\|_{L^\infty_T\dot B^{\frac3q-1}_{q,1}}.\label{equ:Z_T-f}
\end{align}

{\bf Step 4.} The completion of the proof\vspace{0.1cm}

Due to
$\max\Big\{\frac{2}{\widetilde{p}}-\frac1q-\frac1p,\,\frac35\Big(\frac1{q}-\frac1{p}\Big)\Big\}<\epsilon$,
we can take $N$ big enough such that
\begin{align*}\label{}
\frac{T^\f122^{N(\frac3{\widetilde{p}}-\frac3p+1)}}{C(N)}\ll 1 ,\quad
\frac{T2^{N(\frac3q-\frac3p+2)}}{C(N)}\ll 1
\end{align*}
for any $T\le T_0=2^{-2N(1+\epsilon)}.$ Then by a continuous argument, we infer from (\ref{equ:XT-h}), (4.14), (\ref{equ:YT-h-f}) and (\ref{equ:Z_T-f}) that
for any  $T\le T_0$,
\beno
&&X_T\le \f {C_1\big(1+T2^{N(\f3q-\f3p+2)}\big)} {C(N)}, \quad Y_T\le \f {C_2\big(1+T2^{N(\f3q-\f3p+2)}\big)} {C(N)},\\
&&Z_T\le \f {C_3\big(1+T^\f122^{N(\frac3{\widetilde{p}}-\frac3p+1)}+T2^{N(\f3q-\f3p+2)}\big)} {C(N)},
\eeno
where the constants $C_1, C_2, C_3$ are independent of $N$.
Then we deduce from (\ref{equ:theta1final}) and (\ref{equ:u0-h}) that
\begin{align*}
\|\theta(T_0)\|_{\dot B^{\frac3p-2}_{p,1}}&\ge\|\theta_1(T_0)\|_{\dot B^{\frac3p-2}_{p,1}}-\|\Theta_0(T_0)\|_{\dot B^{\frac3p-2}_{p,1}}
-\|\theta_2(T_0)\|_{\dot B^{\frac3p-2}_{p,1}}\nonumber\\&\ge c2^{N(2-\frac3q-\frac3p-3\epsilon)}-
C2^{-\frac{N}{2}(\frac3q-\frac3p+\epsilon)}-C2^{N(\frac3{2q}-\frac3{2p}-\frac52\epsilon)}-C2^{N(\f 3 {\widetilde{p}}-\frac3{2q}-\frac3{2p}-\frac32\epsilon)}
\nonumber\\&\ge\frac c22^{N(2-\frac3q-\frac3p-3\epsilon)}.
\end{align*}
This completes the proof of Theorem \ref{thm:illposed-heat}.

\section*{Acknowledgments}

Qionglei Chen and Changxing Miao  were partially
supported by the NSF of China under grants 11171034 and 11171033.
Zhifei Zhang is partly supported by NSF of China under Grant 10990013 and 11071007,
Program for New Century Excellent Talents in University and Fok Ying Tung Education Foundation.

\end{document}